\documentclass[12pt,reqno]{amsart}
\usepackage{amsfonts}

\newcommand{\Ff}{{\mathbb F}}
\newcommand{\Zz}{{\mathbb Z}}
\newcommand{\Cc}{{\mathbb C}}

\newtheorem{Theorem} {Theorem} [section]
\newtheorem{Proposition} [Theorem] {Proposition}
\newtheorem{Lemma} [Theorem] {Lemma}
\newtheorem{Corollary} [Theorem] {Corollary}
\newtheorem{Conjecture}[Theorem]{Conjecture}
\newtheorem{Example}[Theorem]{Example}

\newcommand{\Proof}{ \noindent{\bf Proof:}\quad }

\def\Tr{\operatorname{Tr}}
\def\tr{\operatorname{tr}}
\def\Im{\operatorname{Im}}

\def\om{\omega}
\def\dmid{\,||\,}
\def\QED{\qed\medskip\par}
\def\AM{\boldsymbol A}

\def\IM{\boldsymbol I}

\catcode`\@=11
\def\sideset#1#2#3{%
  \@mathmeasure\z@\displaystyle{#3}%
  \global\setbox\@ne\vbox to\ht\z@{}\dp\@ne\dp\z@
  \setbox\tw@\box\@ne
  \@mathmeasure4\displaystyle{\copy\tw@#1}%
  \@mathmeasure6\displaystyle{#3{#2}}%
  \dimen@-\wd6 \advance\dimen@\wd4 \advance\dimen@\wd\z@
  \hbox to\dimen@{}\mathop{\kern-\dimen@\box4\box6}%
}
\catcode`\@=12

\numberwithin{equation}{section}

\def\eqref#1{(\ref{#1})}

\def\fl#1{\left\lfloor#1\right\rfloor}

\def\sep{\!\!\mid\!\!}
\def\Punkt{\hbox to5.87494pt{\hfil.\hfil}}
\def\tx{{\tilde x}}
\def\ta{{\tilde a}}
\def\tb{{\tilde b}}
\def\tc{{\tilde c}}
\def\ty{{\tilde y}}
\def\tz{{\tilde z}}

\begin{document}

\title{Gauss Sums, Jacobi Sums, and $p$-ranks of Cyclic Difference Sets}

\author{Ronald Evans, Henk Hollmann, Christian Krattenthaler, Qing Xiang}

\address{
Department of Mathematics,
University of California, San Diego,
La Jolla, California 92093-0112, USA\\
}
\email{\tt revans@euclid.ucsd.edu}

\address{
Philips Research Laboratories, 
Prof.\ Holstlaan 4, 5656 AA Eindhoven,
The Netherlands\\
}
\email{{\tt hollmann@natlab.research.philips.com} }

\address{
Institut f\"ur Mathematik,
Universit\"at Wien,
Strudlhofgasse 4, A-1090 Vienna, Austria\\
}
\email{{\tt kratt@pap.univie.ac.at}\newline\leavevmode\indent
{\it WWW}: {\tt http://radon.mat.univie.ac.at/People/kratt}
}

\address{
Department of Mathematical Sciences,
University of Delaware, Newark, DE 19716, USA\\
}

\email{xiang@math.udel.edu}


\begin{abstract}
We study 
quadratic residue difference sets, GMW difference sets, and
difference sets arising from monomial hyperovals, all of which are
$(2^d-1, 2^{d-1}-1, 2^{d-2}-1)$ cyclic difference sets in
the multiplicative group of the finite field
$\Ff_{2^d}$ of $2^d$ elements, with $d \geq 2.$
We show that, except for a few 
cases with small $d$, these difference sets are 
all pairwise inequivalent.  This is accomplished in part by 
examining their 2-ranks.
The 2-ranks of all of these difference sets were previously known,
except for those connected with the Segre and Glynn hyperovals.
We determine
the 2-ranks of the difference sets arising from the Segre and Glynn
hyperovals, in the following way.  Stickelberger's theorem for
Gauss sums is used to reduce the computation of these 2-ranks to a 
problem of counting certain cyclic binary strings of length $d$.
This counting problem is then solved combinatorially, with the aid of the 
transfer
matrix method.  We give further applications of the 2-rank formulas,
including the determination of the nonzeros of certain binary cyclic codes,
and a criterion in terms of the trace function to decide for which
$\beta$ in  $\Ff_{2^d}^*$ the polynomial $x^6 + x + \beta$ has a zero
in  $\Ff_{2^d}$, when $d$ is odd.
\end{abstract}

\maketitle

\section{Introduction}

Let $G$ be a finite (multiplicative) group of order $v$. A $K$-element
subset $D$ of 
$G$ is called a {\em $(v,K,\lambda)$ difference set\/} in $G$ if the list of
``differences'' $d_1d_2^{-1}$, $d_1,d_2\in D$, $d_1\neq d_2$,
represents each nonidentity element in $G$ exactly $\lambda$
times. Thus, $D$ is a
$(v,K,\lambda)$ difference set in $G$ if and only if it satisfies the
following equation in $\Zz[G]$: 
\begin{equation}\label{eq1.0}
\bigg(\sum _{d\in D} ^{}d\bigg) \bigg(\sum _{d\in D} ^{}d^{-1}\bigg)
=(K-\lambda)1_G+\lambda \sum _{g\in G} ^{}g,
\end{equation}
where $1_G$ is
the identity element of $G$. If the group $G$ is cyclic, then 
$D$ is called a {\em cyclic} difference set.

We say that the $(v,K,\lambda)$ difference sets $D_1$ and $D_2$ 
in an Abelian group $G$ are {\em equivalent\/} if there
exists an automorphism $\alpha$ of $G$ and an element $a\in G$ such that
$\alpha(D_1)=D_2a$. In particular,
if $G$ is cyclic, then $D_1$ and $D_2$ are equivalent if there
exists an integer $t$, $(t,v)=1$, such that $D_1=D_2^{(t)}a$ for some
$a\in G$, where $D_2^{(t)}=\{d^t\mid d\in D_2\}$. 

In the case $G$ is Abelian, using the Fourier inversion formula, we
obtain the following standard result in the theory of difference sets
\cite[p.~323]{turyn}. 

\begin{Lemma}\label{lem1}
Let $G$ be an Abelian group of order $v$. A $K$-subset $D$ is a
$(v,K,\lambda)$ difference set in $G$ if and only if
\begin{equation} \label{eq1.1} 
\chi(D)\overline{\chi(D)}= {K-\lambda}
\end{equation}
for every
nontrivial complex multiplicative character $\chi$ of $G$. Here, $\chi(D)$
stands for $\sum _{d\in D} ^{}\chi(d)$.  
\end{Lemma}

Difference sets are the same objects as symmetric designs with a
regular automorphism group (see \cite[p.~243]{ju}). 
Therefore they play an important role in
the theory of combinatorial designs. The study of difference sets is
also closely related to coding theory, finite geometry, and
communication signal designs. We refer the reader to the book of
Lander \cite{lan} and the paper of Jungnickel \cite{ju} for a survey
of this subject.  

Let $\Ff_q$ be the finite field with $q$ elements, $q$ being a power of the 
prime $p$. Let $\xi_p$ be a fixed complex primitive $p$th root of unity and
let $\Tr _{q/p}$ be the trace from $\Ff_q$ to $\Zz/p\Zz$. Define 
$$\psi\colon\Ff_q\to \Cc^{*},\qquad\psi(x)=\xi_p^{\Tr _{q/p}(x)},$$
which is easily seen to be a nontrivial character of the additive group
of $\Ff_q$. Let
$$\chi\colon\Ff_q^{*}\to \Cc^*$$
be a multiplicative character of $\Ff_q^{*}$. We extend $\chi$ to all
of $\Ff_q$ 
by setting $\chi(0)=0$. Note that $\chi^{q-1}=1$ (the trivial
character), so the order of $\chi$ is prime to $p$. 

Define the {\em Gauss sum} by
$$g(\chi)=\sum_{a\in {\mathbb F}_q}\chi(a)\psi(a).$$
One of the elementary properties of Gauss sums is \cite[Theorem~1.1.4]{be1}
\begin{equation}\label{eq1.2}
g(\chi)\overline{g(\chi)}=q,\hspace{0.1in}{\rm if}\hspace{0.1in}
\chi\neq 1.
\end{equation}

If $\chi_1, \chi_2$ are two multiplicative characters, we define the
{\em Jacobi sum} by
$$J(\chi_1,\chi_2)=\sum_{a\in {\mathbb F}_q}\chi_1(a)\chi_2(1-a).$$
Jacobi sums are closely related to Gauss sums. In fact if $\chi_1\neq
1, \chi_2\neq 1$, and $\chi_1\chi_2\neq 1$, then \cite[Theorem~2.1.3]{be1}
\begin{equation}\label{eq1.3a}
J(\chi_1,\chi_2)=g(\chi_1)g(\chi_2)/g(\chi_1\chi_2),
\end{equation}
and
\begin{equation}\label{eq1.3}
J(\chi_1,\chi_2)\overline{J(\chi_1,\chi_2)}=q.
\end{equation} 

The apparent analogy between \eqref{eq1.1} and $\{$\eqref{eq1.2}, 
\eqref{eq1.3}$\}$ suggests that
there may be a connection between difference sets and
Gauss and Jacobi sums. This is indeed the case. In fact, several authors
have investigated difference sets by using Gauss sums and Jacobi sums;
for example, see \cite{be}, \cite[Chapter~5]{be1}, 
\cite{kYa1}.  However, such investigations deal mostly
with power residue difference sets. In Section~2 of this paper, 
we use Jacobi sums to give
a simple new proof that certain sets which Maschietti
\cite{ma} constructed from hyperovals are cyclic
$(2^d-1,2^{d-1}-1,2^{d-2}-1)$ difference sets; see Theorem~\ref{thm2}. 
Then in Section~3, we use Stickelberger's theorem on the prime ideal
decomposition of Gauss sums to give a new proof of the $p$-rank formula
for Singer difference sets (which include the difference sets arising
from regular and translation hyperovals); see Theorem~\ref{thm6}.
The main result of Section~3 is Theorem~\ref{thm7}, where,
again using Stickelberger's theorem,
we reduce the problem of computing the
2-ranks of the difference sets arising from Segre and Glynn hyperovals
to a problem of counting certain cyclic binary sequences.
In Section~4, we solve the counting problem completely,
with the aid of the transfer matrix method;
see Theorem~\ref{thm:number-solutions}. In the case of the Segre
hyperovals, this yields an explicit formula 
for their 2-ranks in terms of Fibonacci numbers; see 
Corollary~\ref{cor1}. This explicit formula had been
conjectured by Xiang \cite{qx}. In the case of the Glynn hyperovals
we obtain 5-term recurrence relations for the 2-ranks of the
corresponding cyclic difference sets; see 
Theorems~\ref{conj3} and \ref{conj2}.

As applications, in Section~5, we
determine the nonzeros of the binary cyclic codes which result from
the Segre hyperovals, and we give a criterion to decide when the 
equation $x^6+x+\beta=0$, $\beta\in \Ff_{2^d}^*$,
has two distinct solutions in $\Ff_{2^d}$, $d$ odd.

The main applications of our results are contained in Section~6.
There we show that,
except for a few small values of $d$,
the Singer
difference sets and the difference sets which arise from the Segre and
the Glynn hyperovals are all inequivalent; see
Theorem~\ref{thm8}. This answers a question raised by
Maschietti \cite{ma}. We also show that difference
sets arising from hyperovals are inequivalent to quadratic residue
difference sets \cite[p.~244]{ju} (see Theorem~\ref{thm9})
and to
GMW difference sets \cite{go} (see Theorem~\ref{thm10}).

Finally, in Section~7, we use the machinery developed in this paper to
compute the 2-ranks of certain circulant matrices. These matrices are 
closely related to 
incidence matrices corresponding to difference sets arising from hyperovals.

\section{Gauss Sums, Jacobi Sums and Cyclic Difference Sets}

Singer (1938) (cf.~\cite[p.~244]{ju}) 
discovered a large class of difference sets 
which are related to finite projective geometry
and to (generalized) Reed-Muller codes (cf.~\cite[p.~180]{asskey}).
These difference sets $L_0$, whose construction
we give below, have parameters 
$$v=\frac {q^{d}-1} {q-1},\qquad K=\frac {q^{d-1}-1} {q-1},\qquad
\lambda=\frac {q^{d-2}-1} {q-1}$$ 
for $d\geq 2$ and they exist whenever $q$ is a prime power.

The construction is as follows.
Let $\Ff_{q^{d}}$ be the finite field of $q^{d}$ elements, $d\ge2$, and let
$\Tr $ be the trace from $\Ff_{q^{d}}$ 
to $\Ff_q$. We may take a system $L$ of coset representatives of
$\Ff_q^*$ in $\Ff_{q^{d}}^*$
such that $\Tr $ maps $L$ into $\{0,1\}$.
Write $L=L_0\cup L_1$, 
where
$$L_0=\{x\in L\mid \Tr (x)=0\},\qquad L_1=\{x\in L\mid \Tr (x)=1\}.$$

The following proof that $L_0$ is a difference set is due to Yamamoto
\cite{kYa2}. We present it here since it is less well-known than the
standard one using a Singer cycle, and since we will use it
to give a new proof of the $p$-rank formula for the Singer difference
sets (see Theorem~\ref{thm6}). 
This proof of Theorem~\ref{thm1} is essentially an evaluation
of the Eisenstein sum $\chi(L_0)$; see \cite[pp.~389, 400]{be1}.

\begin{Theorem}\label{thm1}
With the above notation, $L_0$ is a $(\frac {q^{d}-1} {q-1},\frac
{q^{d-1}-1} {q-1},\frac {q^{d-2}-1} {q-1})$ difference set in the
quotient group $\Ff_{q^{d}}^*/\Ff_q^*$. 
\end{Theorem}

{\Proof} Let $\chi$ be a nontrivial multiplicative character of
$\Ff_{q^{d}}$ whose restriction to $\Ff_{q}^*$ is
trivial, so that we may view $\chi$ as a character of the quotient group
$\Ff_{q^{d}}^*/\Ff_q^*$. It is easy to see that every nontrivial character of
$\Ff_{q^{d}}^*/\Ff_q^*$ can be obtained in this manner. 
We have
\begin{eqnarray*}
g(\chi)&=&\sum_{y\in
\Ff_{q^{d}}^*}\chi(y)\xi_p^{\Tr _{q^d/p}(y)}=\sum_{a\in
\Ff_{q}^*}\sum_{x\in L}\chi(xa)\xi_p^{\Tr _{q^d/p}(xa)}\\ 
       &=&\sum_{x\in L}\chi(x)\sum_{a\in
\Ff_{q}^*}\chi(a)\xi_p^{\Tr _{q/p}(a\Tr (x))}\\ 
       &=&(q-1)\sum_{x\in L_0}\chi(x)-\sum_{x\in L_1}\chi(x)=
q\sum _{x\in L_0} ^{}\chi(x)=
q\chi(L_0).\\
\end{eqnarray*}
By \eqref{eq1.2}, we have $\chi(L_0)\overline{\chi(L_0)}=q^{d-2}$. The theorem
now follows from Lemma~\ref{lem1}. 
\QED

\noindent{\bf Remark:} When $q=2$, Singer difference sets have
parameters $(2^d-1,2^{d-1}-1,2^{d-2}-1)$, which are a special case of
the {\it Hadamard} parameters $(4n-1,2n-1,n-1)$ (see \cite[p.~244]{ju}).  

Cyclic $(2^d-1,2^{d-1}-1,2^{d-2}-1)$ difference sets have been
extensively studied because of their important applications in
communication signal designs (see \cite{gol}). Known examples of these
difference sets include Singer difference sets, GMW difference sets
(see \cite{go}), quadratic residue difference sets \cite[p.~244]{ju} and 
the cyclic $(2^d-1,2^{d-1}-1,2^{d-2}-1)$ 
difference sets constructed recently by
Maschietti \cite{ma} by using monomial hyperovals. 

We proceed to discuss Maschietti's difference sets.
An $m$-{\it arc} in the projective plane $PG(2,q)$, with $q$ a prime power,
is a set of $m$ points, no three of which are collinear. The maximum
value of $m$ is $q+1$ or $q+2$, according as $q$ is odd or even. If
$q$ is odd, $(q+1)$-arcs are called {\it ovals}. A celebrated theorem
of Segre \cite{segre} (see also \cite[p.~168]{hi}) 
states that all such ovals are given algebraically by
irreducible conics. If $q$ is even, $(q+2)$-arcs are called {\it
hyperovals}. The classic example of a hyperoval is the union of an
irreducible conic and its nucleus (see \cite[p.~164]{hi}). 

Two hyperovals are said to be {\em projectively
equivalent\/} if one hyperoval can be transformed into the other
by a projective linear transformation, i.e., by an element from
$PGL(3,q)$.
By the Fundamental Theorem of Projective Geometry, the group
$PGL(3,q)$ of projective linear transformations of $PG(2,q)$ is transitive on
quadrangles. Thus every hyperoval can be mapped by an element of
$PGL(3,q)$ to a hyperoval containing the {\em fundamental quadrangle}
$(1,0,0), (0,1,0), (0,0,1)$, and $(1,1,1)$. 
{}From now on, we will restrict our attention to those hyperovals
$PG(2,q)$ with $q>2$ which contain the fundamental quadrangle. The
following result of Segre (see \cite[Theorem~8.4.2]{hi}) 
shows that any such hyperoval can be expressed in terms of a
{\em permutation polynomial} over $\Ff_q$, i.e., a polynomial which, when
interpreted as a function, permutes $\Ff_q$.

\begin{Theorem}{\em (Segre).} \label{thm3'}
Let $q>2$ be a power of $2$.
Then any hyperoval in
$PG(2,q)$ containing the fundamental quadrangle 
is equal to a $(q+2)$-arc
$$D(f)=\{(1,t,f(t))\mid t\in \Ff_q\}\cup \{(0,1,0), (0,0,1)\},$$ 
where $f$ is a
permutation polynomial over $\Ff_q$ of degree at most $q-2$,
satisfying $f(0)=0, f(1)=1$, and such that for each $s\in \Ff_q$, 
$$f_s(x)=\begin{cases} \dfrac {f(x+s)+f(s)} {x},&\text{if }x\ne0,\\
0,&\text{if }x=0,\end{cases}$$
is a permutation polynomial. 

Conversely, every such set $D(f)$ is a hyperoval. 
\end{Theorem}

These polynomials $f(x)$ are also called
$o$-{\it polynomials}. 

\medskip
A subclass of hyperovals is the set of {\it monomial hyperovals}. These are
hyperovals in $PG(2,2^d)$ projectively equivalent to
$$D(x^k)=\{(1,t,t^k)\mid t\in \Ff_{2^d}\}\cup \{(0,1,0), (0,0,1)\}.$$ 
Known examples of monomial hyperovals in $PG(2,2^d)$ include 
{\leftskip10pt\parindent0pt

\medskip
({\refstepcounter{equation}\label{reg-oval}}\ref{reg-oval}) 
the {\it regular} hyperoval $D(x^2)$, 

({\refstepcounter{equation}\label{trans-oval}}\ref{trans-oval}) 
the {\it translation} hyperovals $D(x^{2^i})$, where
$(d,i)=1$, $1<i<d/2$,

({\refstepcounter{equation}\label{Segre-oval}}\ref{Segre-oval}) 
the {\it Segre} hyperoval $D(x^6)$, where $d\geq 5$ is odd, 

\hangindent29pt\hangafter1
({\refstepcounter{equation}\label{Glynn-I-oval}}\ref{Glynn-I-oval}) 
the {\it Glynn} type (I) hyperovals $D(x^{\sigma+\gamma})$,
where $d\geq 7$ is odd, $\sigma=2^{ (d+1)/ {2}}$, and
$\gamma=2^{(3d+1)/4}$ if $d\equiv1$ mod 4, whereas $\gamma=2^{(d+1)/4}$ 
if $d\equiv3$ mod 4,

({\refstepcounter{equation}\label{Glynn-II-oval}}\ref{Glynn-II-oval}) 
the {\it Glynn} type (II) hyperovals $D(x^{3\sigma +4})$,
where $d\geq 11$ is odd, and $\sigma=2^{ (d+1)/ {2}}$. 
\medskip\par}

\noindent{\bf Remark:} 
We have not listed the translation hyperovals for $i=1$ nor the Segre
and Glynn hyperovals for $d=3$, as these are all projectively equivalent 
to the regular hyperovals. 
The translation hyperoval for a given $i>d/2$ is equivalent to the one
where $i$ is replaced by $d-i$. The Glynn hyperovals for $d=5$ are
equivalent to translation hyperovals, while the Glynn type (II)
hyperovals for $d=7$ and $d=9$ are respectively equivalent to the
Segre hyperoval $D(x^6)$ and the Glynn type (I) hyperoval
$D(x^{160})$. These equivalences all follow from Lemma~\ref{lem3'} below.
It can be shown using Lemma~\ref{lem3'} that the hyperovals listed in
\eqref{reg-oval}--\eqref{Glynn-II-oval} are all mutually projectively 
inequivalent.  For example,
to show that the Segre hyperovals are inequivalent to the translation
hyperovals, one shows that none of $6$, $1/6$, $-5$, $-1/5$, $6/5$, or
$5/6$
(mod $q-1$) can equal a power of 2 (mod $q-1$), where $q=2^d$ for odd
$d \geq 5$.
The computations involved are tedious but straightforward.

Glynn \cite{gl} conjectured that the hyperovals listed in 
\eqref{reg-oval}--\eqref{Glynn-II-oval}
(together with the equivalent hyperovals given by Lemma~\ref{lem3'} below)
comprise {\em all\/} monomial hyperovals. 
This conjecture still remains open. Some progress was made
in a recent paper by Cherowitzo and Storme \cite{cs}. 

\begin{Lemma}\label{lem3'}{\em (\cite{gl}, \cite{pt})}
Let $q>2$ be a power of $2$.
Two monomial hyperovals $D(x^j)$ and $D(x^k)$ in $PG(2,q)$ are
projectively equivalent if and only if 
$j\equiv k$, $1/k$, $1-k$, $1/(1-k)$, $k/(k-1)$, or $(k-1)/k$ {\em(mod $q-1$)}.
\end{Lemma}
{\Proof} 
First, if $D(x^k)$ is a monomial hyperoval in $PG(2,q)$,
then it is well known that
for any $j$ congruent to $k$, $1/k$, $1-k$, $1/(1-k)$, 
$k/(k-1)$, or $(k-1)/k$ modulo $q-1$, $D(x^j)$ is projectively equivalent to 
$D(x^k)$ (see \cite{gl}). 
Conversely, 
it follows from \cite[Theorem~12.5.3]{pt} that
if $D(x^j)$ is projectively equivalent to $D(x^k)$, then $j$ is 
congruent to one of $k$, $1/k$, $1-k$, $1/(1-k)$, $k/(k-1)$, $(k-1)/k$ modulo 
$q-1$.   
This completes the proof.
\QED

Let $\tau\colon \Ff_{2^d}\to \Ff_{2^d}$ be defined by 
$$\tau(x)=x+x^k,$$  
and let
$\Im\tau$ be the image of the map $\tau$. We need the following lemma
from Maschietti \cite{ma}. For the convenience of the reader, we
include a proof. 
In the statement of the lemma, and its proof, we use the convention
that points in $PG(2,q)$ are labelled as $(z,x,y)$.

\begin{Lemma}\label{lem2}
Let $q=2^d$.
The $(q+2)$-set $D(x^k)=\{(1,t,t^k)\mid t\in \Ff_q\}\cup \{(0,1,0),
(0,0,1)\}$ in $PG(2,q)$ is a hyperoval if and only if
$(k,q-1)=1$ and $\tau$ is a two-to-one map from $\Ff_q$ into itself. 
\end{Lemma}

{\Proof} If $D(x^k)$ is a hyperoval, then $x^k$ is a
permutation polynomial on $\Ff_q$, so $(k,q-1)=1$. Furthermore,
any affine line $y=x+a$ ($a\in \Ff_q$) meets $D(x^k)$ at either 0 or 2
distinct points, which implies that $x^k+x+a=0$ has either 0 or 2 distinct
solutions in $\Ff_q$, that is, the map $\tau$ is two-to-one. 

Conversely, if $\tau$ is two-to-one, then the equation $x^k+x=0$
has exactly two solutions, namely 0 and 1. Hence, there is no $(k-1)$th
root of 1 distinct from 1 in $\Ff_q$, so that $(k-1, q-1)=1$.  

It is straightforward to show that a line 
of $PG(2,q)$ with homogeneous equation $cy=ax+bz$, with $abc=0$, 
intersects $D(x^k)$ at either 0 or 2 points. It remains to show that
the line $y=ax+bz$, with $ab\ne0$,
intersects $D(x^k)$ at either 0 or 2 points. 

It is clear that there are no intersection points $(0,x,y)$ of the
line with $D(x^k)$. Thus set $z=1$ and consider the affine equation
$x^k=ax+b$. We must show that this equation has 0 or 2 solutions. This follows
from the hypothesis that $\tau$ is two-to-one, as can be seen after
making the change of variable
$x=a^{\frac {1} {k-1}}x'$. This completes the proof.
\QED

\noindent{\bf Remark:} 
{}From Lemma~\ref{lem2} and its proof, it follows that when $D(x^k)$ is a
hyperoval in $PG(2,q)$, with $q=2^d$, then $k$ and $k-1$ are both
relatively prime to $q-1$.

\medskip
Using the geometry of hyperovals, Maschietti \cite{ma} proved that if
$D(x^k)$ is a hyperoval in $PG(2,2^d)$, then 
$$D_{k,d}=\Im\tau\setminus\{0\}$$ 
is a $(2^d-1,2^{d-1}-1,2^{d-2}-1)$ difference
set in $\Ff_{2^d}^*$. In the following theorem, we give a simpler
proof using Jacobi sums.  

\begin{Theorem}\label{thm2}
Let $q=2^d$.
If $D(x^k)=\{(1,t,t^k)\mid t\in \Ff_{q}\}\cup \{(0,1,0), (0,0,1)\}$ is
a hyperoval in $PG(2,q)$, then 
$D_{k,d}=\Im\tau\setminus \{0\}$ is a $(q-1,q/2-1,q/4-1)$ difference set in
$\Ff_q^*$. 
\end{Theorem}

{\Proof} Let $\chi$ be a nontrivial multiplicative character of
$\Ff_q$. Since $D(x^k)$ is a hyperoval, $\tau$ is
two-to-one by Lemma~\ref{lem2}, and we have 
$$\chi(D_{k,d})=\frac {1} {2}\sum_{x\in \Ff_q}\chi(x+x^k)=\frac {1}
{2}\sum_{x\in \Ff_q}\chi(x)\chi(1+x^{k-1}).$$ 
By the previous remark, $(k-1,q-1)=1$, so there exists a
multiplicative character $\phi$ of $\Ff_q$ such that
$\chi=\phi^{k-1}$. Hence 
\begin{equation}\label{eq2.1}
\chi(D_{k,d})=\frac {1} {2}\sum_{x\in
\Ff_q}\phi(x^{k-1})\chi(1+x^{k-1})=\frac {1} {2} J(\phi,\chi).
\end{equation}
Noting that $\phi\chi=\phi^k$ is nontrivial, we have, by \eqref{eq1.3}, 
$\chi(D_{k,d})\overline{\chi(D_{k,d})}=2^{d-2}$. Now the theorem
follows from Lemma~\ref{lem1}. 
\QED

The following proposition shows that 
 two projectively equivalent monomial hyperovals give rise to two
 equivalent cyclic difference sets under the construction of
 Theorem~\ref{thm2}. 

\begin{Proposition}
Let $q=2^d$, $q>2$.
If $D(x^k)$ and $D(x^j)$ are two projectively equivalent monomial
hyperovals in $PG(2,q)$,  then the corresponding difference sets
$D_{k,d}$ and $D_{j,d}$ constructed in
Theorem~{\em\ref{thm2}} are equivalent. 
\end{Proposition}

{\Proof} By Lemma~\ref{lem2}, we know that $j$ is congruent to one of
$k$, $1/k$, $1-k$, $1/(1-k)$, $k/(k-1)$, $(k-1)/k$ modulo
$q-1$. Now direct calculation shows that 
\begin{gather}\label{eq:first}
D_{\frac {k} {k-1},d}=D_{\frac {k-1} {k},d}=D_{k,d}^{(k-1)},\\
D_{1-k,d}=D_{\frac {1} {1-k},d}=D_{k,d}^{(\frac {1-k} {k})},\\
D_{\frac {1} {k},d}=D_{k,d},
\end{gather}
where $D^{(t)}=\{d^t\mid d\in D\}$. For example, to prove
\eqref{eq:first}, write $c=x+x^k$, make the variable change $x=c/y$,
then multiply both sides by $y^k$.
This completes the proof of the proposition.
\QED    

The difference sets $D(x^{2^i})$ arising from the regular and translation
hyperovals are the Singer difference sets 
$\{y\in \Ff_{2^d}^* \mid \Tr(y)=0\}$.
This follows from Lemma~\ref{lem2} and Hilbert's Theorem~90 
\cite[p.~56]{lidl}.  It is an interesting
problem to determine whether the cyclic difference sets arising from
the Segre and Glynn hyperovals are inequivalent to each other and 
to the previously
known ones (such as Singer difference sets, quadratic residue
difference sets, GMW difference sets). 
We approach
these problems by studying the 2-ranks (as defined in Section~3)
of these cyclic difference sets. Our conclusions are presented in
Theorems~\ref{thm8}, \ref{thm9}, and \ref{thm10}.

\section{Stickelberger's Theorem and the $p$-ranks of Cyclic Difference Sets}

We start this section by reviewing a few standard notions from design
theory. We refer the reader to \cite[p.~243]{ju} for details and
more information.

Let $G$ be a (multiplicative) Abelian group of order $v$, and let $D$ be a
$(v,K,\lambda)$ difference set in $G$. Then ${\mathcal D}=({\mathcal
P}, {\mathcal B})$ is a $(v,K,\lambda)$ symmetric design with a
regular automorphism group $G$, where the set ${\mathcal P}$ of {\em
points} of $\mathcal D$ is $G$,
and where the set ${\mathcal B}$ of {\em blocks} of $\mathcal D$ is
$\{gD\mid g\in G\}$. (As it turns out, the blocks are all distinct.) 
This design is
usually called the {\em development\/} of $D$. The incidence matrix of
${\mathcal D}$ is the matrix $A$ whose rows are indexed by the blocks
$B$ of ${\mathcal D}$ and whose columns are indexed by the points $g$ of
${\mathcal D}$, where the entry $A_{B,g}$ in row $B$ and column $g$ is
1 if $g\in B$, and 0 otherwise. 

The {\em $p$-ary code} of $D$, denoted ${\mathcal C}_p(D)$, is defined to be
the row space of $A$ over $\Ff_p$, the field of $p$ elements. This
code is also the $p$-ary code of ${\mathcal D}$, denoted by ${\mathcal
C}_p({\mathcal D})$. The $\Ff_p$-dimension of ${\mathcal C}_p(D)$ is
usually called the $p$-{\it rank} of the difference set $D$. We shall
also denote it by ${\rm rank}_p{\mathcal D}$. It is
well-known that ${\mathcal C}_p(D)$ is of interest only if
$p\mid (K-\lambda)$; see \cite{BrHaHa} or \cite[p.~297]{ju}. 
So from now on, we always assume that
$p\mid (K-\lambda)$.

The $p$-ranks of difference sets have been studied extensively for
several reasons. First, if $D_1$ and $D_2$ are two equivalent
$(v,K,\lambda)$ difference sets in an Abelian group $G$, 
then the $p$-ranks of $D_1$ and $D_2$ are the same. 
Therefore, the $p$-ranks can help us to distinguish two
inequivalent difference sets.  

Secondly, let $D$ be a cyclic $(v,K,\lambda)$ difference set in
$\Zz/v\Zz$. Corresponding to $D$ is the {\em characteristic
sequence}
$\{a_i\}_{0\leq i\leq v-1}$ given by $a_i=1$ if $i\in D$, and $a_i=0$
otherwise. This sequence has a two-level autocorrelation function 
(see \cite[p.~62]{gol}). The
{\it linear complexity} (or {\em linear span}) of $\{a_i\}$ is the smallest
degree of a linear shift register over $\Ff_2$ which is capable of
generating $\{a_i\}$. It turns out that the linear complexity of
$\{a_i\}$ is the same as the 2-rank of $D$. Therefore, the study of
2-ranks of cyclic difference sets is important for communications and
cryptographic applications; see \cite{gol}.  

A third reason is provided by Theorem~\ref{thm3} below.
MacWilliams and Mann \cite{mm}, Goethals and Delsarte \cite{gd},
and Smith \cite{sm} proved that the $p$-rank of the
$((q^{d}-1)/(q-1),(q^{d-1}-1)/(q-1),(q^{d-2}-1)/(q-1))$ Singer
difference set is ${\binom {p+d-2 }{ d-1}}^s+1$, where $q=p^s$ for
prime $p$. 
Subsequently, Hamada \cite{h} made the following conjecture. 

\begin{Conjecture} 
Let ${\mathcal D}$ be a symmetric design with Singer parameters
$((q^{d}-1)/(q-1),(q^{d-1}-1)/(q-1),(q^{d-2}-1)/(q-1))$, where
$q=p^s$. Then one has 
$${\rm rank}_p{\mathcal D}\geq {\binom {p+d-2 }{ d-1}}^s+1,$$ 
with equality if and only if ${\mathcal D}$ is the development of a classical
Singer difference set. 
\end{Conjecture}

This conjecture still remains open. But Hamada and Ohmori \cite{ho}
proved the following interesting result in this direction.

\begin{Theorem}\label{thm3}
Let ${\mathcal D}$ be a symmetric design with parameters
$(2^d-1,2^{d-1}-1,2^{d-2}-1)$. Then ${\mathcal D}$ is the development
of the cyclic Singer difference set with these parameters if and only
if ${\rm rank}_2{\mathcal D}=d+1$. 
\end{Theorem}

In Theorem~\ref{thm6}, we will give a simple new proof of the $p$-rank
formula for the classical Singer difference sets. 
Our approach is different from previous ones in that it
uses Gauss and Jacobi sums, and Stickelberger's theorem on the prime ideal
factorization of Gauss sums. 

We first quote a result (Theorem~\ref{thm4}) of MacWilliams and Mann
\cite{mm}. For 
the convenience of the reader, we include a proof here. 

Let $B=\sum_{g\in G}b_gg$ be an element of the group algebra $F[G]$,
where $G$ is an Abelian group and $F$ is a field. 
We associate with $B$ the matrix
$(b_{g(g')^{-1}})$, whose rows and columns are labeled by the group
elements $g$ and $g'$. 
The rank of $(b_{g(g')^{-1}})$ is called the {\it rank} of
$B$. We remark that in particular, if $B=\sum _{g\in D} ^{}g$, where
$D$ is a
$(v,K,\lambda)$ difference set in $G$, then the $p$-rank of $D$, i.e.,
the $\Ff_p$-dimension of the code ${\mathcal C}_p(D)$, is the same as the rank
of $B$ over $\Ff_p$. The reason is that the rows and columns of the
incidence matrix of the development of $D$ can be arranged so that the
incidence matrix is the transpose $(b_{g(g')^{-1}})^T$ of the matrix
$(b_{g(g')^{-1}})$.

\begin{Theorem}\label{thm4}
Let $G$ be an Abelian group of order $v$ and exponent $v^*$, let $F$
be a field of characteristic $p$ not dividing $v$ which contains the
$v^*$th roots of unity, and let $B=\sum_{g\in G}b_gg$ be an element of
the group algebra $F[G]$. Then the rank of $B$ over $F$ is equal to
the number of characters $\chi :G\rightarrow F^*$ satisfying
$\chi(B)\neq 0$, where $\chi(B)=\sum _{g\in G} ^{}b_g\chi(g)$.
\end{Theorem}

{\Proof} Let $\big(\chi(g)\big)$ be a matrix whose rows are labeled by the $v$
characters $\chi$ and whose columns are labeled by the $v$ group
elements $g$, so that the entry in row $\chi$ and column $g$ is
$\chi(g)$. 
This matrix is nonsingular since $\frac {1}
{v}\big(\chi(g)\big)\big(\chi(g^{-1})\big)^{T}$ is the identity matrix.
The result now follows from the identity
$$\big(\chi(g)\big)\big(b_{g(g')^{-1}}\big)\big(\chi(g^{-1})\big)^{T}=
v\cdot {\rm diag}\big(\chi(B)\big).$$
\QED 

We will also need Stickelberger's result (Theorem~\ref{thm5}) 
on the prime ideal
decomposition of Gauss sums. We first introduce some notation. 

Let $p$ be a prime, $q=p^s$, and let $\xi_{q-1}$ be a (complex) primitive
$(q-1)$th
root of unity. Fix any prime ideal $\mathfrak{p}$ in $\Zz[\xi_{q-1}]$ lying
over $p$. Then $\Zz[\xi_{q-1}]/\mathfrak{p}$ is a finite field of order
$q$, which we identify with $\Ff_q$. Let $\omega_{\mathfrak{p}}$
be the Teichm\"uller character on $\Ff_q$, i.e., an isomorphism
$$\omega_{\mathfrak{p}}: \Ff_q^{*}\rightarrow
\{1,\xi_{q-1},\xi_{q-1}^2,\dots ,\xi_{q-1}^{q-2}\}$$ 
satisfying 
\begin{equation}\label{eq3.1}
\omega_{\mathfrak{p}}(\alpha)\quad ({\rm
mod}\hspace{0.1in}{\mathfrak{p}})=\alpha,
\end{equation}  
for all $\alpha$ in $\Ff_q^*$.
The Teichm\"uller character $\omega_{\mathfrak{p}}$ has order $q-1$;
hence it generates all multiplicative characters of $\Ff_q$. 

Let $\mathfrak{P}$ be the prime ideal of $\Zz[\xi_{q-1},\xi_p]$ lying above
$\mathfrak{p}$. For an integer $a$, let
$s(a)=v_{\mathfrak{P}}(g(\omega_{\mathfrak p}^{-a}))$, where $v_{\mathfrak{P}}$ 
is
the $\mathfrak{P}$-adic valuation. Thus ${\mathfrak
P}^{s(a)}\dmid g(\omega_{\mathfrak p}^{-a})$. The following evaluation of
$s(a)$ is due to Stickelberger (see \cite[p.~344]{be1}, \cite[p.~96]{wa}).

\begin{Theorem}\label{thm5}
Let $p$ be a prime, and $q=p^s$. 
For an integer $a$ not divisible by $q-1$, 
let $a_0+a_1p+a_2p^2+\cdots +a_{s-1}p^{s-1}$, $0\leq a_i\leq p-1$,
be the $p$-adic expansion of the reduction of $a$ modulo $q-1$.
Then
$$s(a)=a_0+a_1+\cdots +a_{s-1},$$
that is, $s(a)$ is the
sum of the $p$-adic digits of the reduction of $a$ modulo $q-1$. 
\end{Theorem}

An easily proved consequence of Theorem~\ref{thm5} is the formula
\cite[p.~347]{be1}, \cite[p.~98]{wa} 
\begin{equation} \label{eq3.1'}
s(a)=\frac {p-1} {q-1}\sum _{i=0} ^{s-1}L(ap^i), 
\end{equation}
where $L(x)$ is the reduction of $x$ modulo $q-1$.

We are now ready to give a new proof of the $p$-rank formula for 
Singer difference sets. 

\begin{Theorem}\label{thm6}
For $q=p^s$,
let $L_0$ be the $(\frac {q^{d}-1} {q-1},\frac {q^{d-1}-1} {q-1},\frac
{q^{d-2}-1} {q-1})$ difference set in the quotient group
$\Ff_{q^{d}}^*/\Ff_q^*$, as in Theorem~{\em\ref{thm1}}. Then the
$p$-rank of $L_0$ is 
$${\binom {p+d-2 }{ d-1}}^s+1.$$
\end{Theorem}

{\Proof} 
For $\mathfrak{p}$ a prime ideal in $\Zz[\xi_{q^d-1}]$ lying over
$p$,
let $\omega_{\mathfrak p}$ be the Teichm\"uller character on $\Ff_{q^d}$
and let $\chi=\omega_{\mathfrak p}^{-(q-1)}$. Then $\chi$ is a generator of the
character group of 
$\Ff_{q^d}^*/\Ff_q^*$. From the proof of Theorem~\ref{thm1}, we know that for
each $a$, $0<a<\frac {q^{d}-1} {q-1}$, 
\begin{equation}\label{eq3.2}
q\cdot \chi^{a}(L_0)=g(\chi^{a}).
\end{equation}
Note that if $\chi_0$ is the trivial character, then
$\chi_0(L_0)=|L_0|\not\equiv 0$ (mod $p$).
Thus by the definition~\eqref{eq3.1} 
of the Teichm\"uller character 
and Theorem~\ref{thm4}, 
the $p$-rank of $L_0$ is $1+A(q,d)$, where $A(q,d)$ is
the number of $\chi^a$, $0<a< (q^{d}-1)/(q-1)$, such
that $\chi^a(L_0) ({\rm mod}\hspace{0.1in} \mathfrak{p})\neq 0$.

Let $\mathfrak{P}$ be the prime of $\Zz[\xi_{q^d-1},\xi_p]$ lying
above $\mathfrak{p}$. 
Since $\mathfrak P\mid\chi^a(L_0)$ if and only if $\mathfrak p\mid
\chi^a(L_0)$, 
$A(q,d)$ is equal to the number of $\chi^a$,
$0<a<(q^{d}-1)/(q-1)$, such that $\mathfrak P
\hbox{${}\not\kern2.5pt\mid{}$}\chi^a(L_0)$.

By the definition of $s(a)$ (with $q^d$ in place of $q$),
$\mathfrak{P}^{s((q-1)a)}\dmid g(\chi^a)$. Since
$\mathfrak{P}^{(p-1)s}\dmid q$, we see from \eqref{eq3.2} that $A(q,d)$
is equal to the number of
$a$, $0<a<(q^{d}-1)/(q-1)$, such that $s((q-1)a)=(p-1)s$, which, in
turn, is equal to the number of $x$, $0<x<q^d-1$, $(q-1)\mid x$, such that
$s(x)=(p-1)s$.  

For $x$ with $0<x<q^d-1$, write 
$x=\sum_{j=0}^{d-1}\sum_{i=0}^{s-1}x_{i,j}p^iq^j$, $0\leq
x_{i,j}<p$. 
By \eqref{eq3.1'} (with $q^d$ in place of $q$), 
\begin{equation} \label{eq3.1''}
s(x)=\frac {p-1} {q^d-1}\sum _{j=0} ^{d-1}\sum _{i=0} ^{s-1}L(xp^iq^j), 
\end{equation}
where $L(y)$ denotes the reduction of $y$ modulo $q^d-1$. Suppose that
$(q-1)\mid x$ and $s(x)= (p-1)s$. Then, by \eqref{eq3.1''},
\begin{equation} \label{eq3.1'''}
(q^d-1)s= \sum _{j=0} ^{d-1}\sum _{i=0} ^{s-1}L(xp^iq^j).
\end{equation}
Since 
$$0<\sum _{j=0} ^{d-1}L(xp^iq^j)\equiv \sum _{j=0} ^{d-1}xp^iq^j\equiv
0\pmod{q^d-1},   $$
it follows that 
$\sum _{j=0} ^{d-1}L(xp^iq^j)\ge q^d-1 $
for each $i$. Thus, by \eqref{eq3.1'''},
\begin{equation}\label{eq3.1''''}
\sum _{j=0} ^{d-1}L(xp^iq^j)= q^d-1 
\end{equation}
for each $i$. Write $x=a_0+a_1q+\dots+a_{d-1}q^{d-1}$, where $a_j=\sum
_{i=0} ^{s-1}x_{i,j}p^i$, $0\le j\le d-1$. 
Then by \eqref{eq3.1''''} with $i=0$, we have
$a_0+a_1+\dots+a_{d-1}=q-1$. Thus $\sum _{j=0} ^{d-1}x_{0,j}\ge
p-1$. Similarly, by \eqref{eq3.1''''} for general $i$, we obtain
\begin{equation} \label{eq3.1'''''}
\sum _{j=0} ^{d-1}x_{i,j}\ge p-1 
\end{equation}
for {\em all\/} $i$. As 
$$(p-1)s\le \sum _{i=0} ^{s-1}\sum _{j=0} ^{d-1}x_{i,j}=s(x)=(p-1)s,$$
it follows that equality holds in \eqref{eq3.1'''''} for each $i$.
The equation $\sum_{j=0}^{d-1}x_{i,j}=p-1$, $0\leq x_{i,j}<p$,
has $\binom {p+d-2}{ d-1}$ solutions for each $i$, $0\le i< s$.
Therefore $A(q,d)={\binom {p+d-2}{
d-1}}^s$. This completes the proof. 
\QED

Next we show how to compute the 2-ranks of those cyclic difference
sets $D_{k,d}$
arising from monomial hyperovals in Section~2. For succinctness,
we will actually state the theorem in
terms of the 2-rank of the complement $\overline{D_{k,d}}$. 
(It follows easily from \eqref{eq1.0} that if $D$ is a $(v,K,\lambda)$
difference set in $G$ then the complement $\overline D$ of $D$ in $G$
is a $(v,v-K,\lambda-2K+v)$ difference set in $G$.)
This is just a matter of convenience, as
the 2-rank of the difference set $D_{k,d}$ is exactly 1 more than the
2-rank of the complement $\overline {D_{k,d}}$ (see the proof of
Theorem~\ref{thm7}).

\begin{Theorem}\label{thm7}
Let $D_{k,d}$ be the $(2^d-1,2^{d-1}-1,2^{d-2}-1)$ cyclic difference
set in $\Ff_{2^d}^*$ constructed from the hyperoval $D(x^k)$ as in
Theorem~{\em \ref{thm2}}, and let $\overline{D_{k,d}}$ be the complement of
$D_{k,d}$ in $\Ff_{2^d}^*$. Then the {\em 2}-rank of $\overline{D_{k,d}}$ is
equal to the number of $a$'s, $0<a<2^d-1$, such that
\begin{equation}\label{eq:keyeq1}
s(a)+s((k-1)a)=s(ka)+1,
\end{equation}
where $s(a)$ is as defined above Theorem~{\em \ref{thm5}} with $q=2^d$.
\end{Theorem}

{\Proof} 
For $\mathfrak{p}$ a prime ideal in $\Zz[\xi_{2^d-1}]$ lying over
$2$, let $\omega=\omega_{\mathfrak p}$ be the Teichm\"uller character on
$\Ff_{2^d}$. By \eqref{eq2.1}, we see that for each $a$,
$0<a<2^d-1$, 
\begin{equation}\label{eq3.3}
-2\cdot
\omega^{-(k-1)a}\left(\overline{D_{k,d}}\right)=J(\omega^{-a},
\omega^{-(k-1)a}).
\end{equation}  
By the definition of $s(a)$ (with $q=2^d$) and \eqref{eq1.3a},  
\begin{equation} \label{eq3.3a} 
\mathfrak{p}^{s(a)+s((k-1)a)-s(ka)}\dmid J(\omega^{-a},\omega^{-(k-1)a}).
\end{equation}
Since $\mathfrak p\dmid 2$, we see from \eqref{eq3.3} that
the number of $a$'s, $0<a<2^d-1$, such that
$\omega^{-(k-1)a}\left(\overline{D_{k,d}}\right)$ (mod $\mathfrak{p}$)
$\neq 0$, is equal to the number of $a$'s, $0<a<2^d-1$, such that
$s(a)+s((k-1)a)=s(ka)+1$. Since the cardinality of
$\overline{D_{k,d}}$ is $2^{d-1}\equiv 0$ (mod 2),
$\chi_0(\overline{D_{k,d}})$ (mod $\mathfrak{p}$) equals 0 
for the trivial character
$\chi_0$. (In contrast, $\chi_0(D_{k,d})$ (mod $\mathfrak{p}$) equals $1$, 
which is why we prefer to
state the theorem in terms of $\overline{D_{k,d}}$.) 
Thus, by Theorem~\ref{thm4}, we
see that the 2-rank of $\overline{D_{k,d}}$ is equal to the number of
$a$'s, $0<a<2^d-1$, for which $s(a)+s((k-1)a)=s(ka)+1$. This
completes the proof. 
\QED 
 
Let $B_k(d)$ be the number of $a$'s, $0<a<2^d-1$, for which
$s(a)+s((k-1)a)=s(ka)+1$. By the above theorem, in order to
compute the 2-rank of ${D_{k,d}}$, we need to compute
$B_k(d)$. We will discuss the determination of $B_k(d)$ 
in the next section. But first of
all, we observe that $B_k(d)$ is always a multiple of $d$. We state
this as a lemma. 

\begin{Lemma}\label{lem3}
Let $B_k(d)$ be the number of $a$'s, $0<a<2^d-1$, for which
$s(a)+s((k-1)a)=s(ka)+1$. Then $d\mid B_k(d)$. 
\end{Lemma}

{\Proof} For any $a$, $0<a<2^d-1$, it follows easily from
\eqref{eq3.1'} that
\begin{equation}\label{eq3.4}
s(a)+s((k-1)a)-s(ka)=\sum_{i=0}^{d-1}\left(\fl{\frac
{2^ika} {2^d-1}}-\fl{\frac {2^ia}
{2^d-1}}-\fl{\frac {2^i(k-1)a} {2^d-1}}\right),
\end{equation} 
where
$\fl x$ denotes the greatest integer $\le x$ (see 
\cite[p.~349]{be1}, \cite[p.~98]{wa}).

Clearly, if $a$ is a solution to $s(a)+s((k-1)a)=s(ka)+1$, then
$a\cdot 2^j$ is also a solution, for all $j$, $0\leq j\leq d-1$. We
contend that $a\cdot 2^j$, $j=0,1,\dots ,d-1$, are necessarily {\it
distinct} modulo $2^d-1$. If they were not distinct, then the sum
on $i$ in \eqref{eq3.4} would have more than one term equal to 1, contradicting
the fact that
$s(a)+s((k-1)a)-s(ka)=1$. This completes the proof. 
\QED

\section{Enumeration of Cyclic Binary Strings and Computation of 2-ranks}

Let $d$ be an integer $\ge 2$.
In this section, we are concerned with determining explicitly the
2-ranks of the difference sets $\overline {D_{k,d}}$ 
corresponding to the hyperovals
$D(x^k)$ (which, by the remark preceding Theorem~\ref{thm7}, is
equivalent to determining the 2-ranks of the difference sets 
$D_{k,d}$ themselves). 
By Theorem~\ref{thm7}, the 2-rank of $\overline{D_{k,d}}$ equals the
number $B_k(d)$ of solutions $a$, $0<a<2^d-1$, to equation
\eqref{eq:keyeq1}. The case of the difference sets arising from the translation
and the regular hyperovals is already covered by Theorem~\ref{thm6}, as
these difference sets are instances of Singer difference sets (see the
last paragraph of Section~2). Hence, we concentrate on the
difference sets arising from the Segre and the Glynn hyperovals. Most
of the time we will not actually deal with the number $B_k(d)$ of
solutions itself, but with the quantity $A_k(d)=B_k(d)/d$. For the
convenience of the reader we have listed the first few values for
$A_6(d)$ (the Segre case, see \eqref{Segre-oval}),
for $A_{\sigma+\gamma}(d)$ (the Glynn type (I) case, see
\eqref{Glynn-I-oval}), and for $A_{3\sigma+4}(d)$ 
(the Glynn type (II) case, see \eqref{Glynn-II-oval}).
These values were originally found by the computer
algebra package MAGMA \cite{cp}. 

\begin{table}[h]
\begin{center}
\begin{tabular}{|c||c|c|c|c|c|c|c|c|c|c|c|c|c|}
\hline
{\sc $d$}& 3 & 5 & 7 & 9 & 11 & 13 & 15 & 17 & 19 & 21 & 23 & 25 \\ \hline

{\sc $A_{6}(d)$} & 1 & 3 & 5 & 9 & 15 & 25 & 41 & 67 & 109 & 177 & 287 &
465\\ \hline 

{\sc $A_{\sigma+\gamma}(d)$} & 1 & 1 & 3 & 7 & 13 & 23 & 45 & 87 & 167 & 321 &
619 & 1193 \\ \hline 

{\sc $A_{3\sigma+4}(d)$} & 1 & 1 & 5 & 7 & 21 & 37 & 89 & 173 & 383 & 777 &
1665 &3441\\ \hline

\end{tabular}
\end{center}
\vskip6pt
\caption{}
\end{table}

When we entered the values from the
table into the {\sl Maple} package {\tt gfun} \cite{SaZi}, we
were able to guess a
linear recurrence for the numbers $A_k(d)$ in each of the
three cases. These recurrences are established 
by rather elaborate combinatorial arguments in
Theorems~\ref{thm:number-solutions}, \ref{conj3} and \ref{conj2}.
Consequently, we obtain explicit
formulas for the 2-ranks of the difference sets $D_{k,d}$ corresponding
to the Segre and Glynn hyperovals $D(x^k)$; see
Corollaries~\ref{cor1}, \ref{cor3} and \ref{cor2}.

Table~1 suggests that for $d\ge15$ the 2-ranks $A_6(d)$ 
for the Segre case are the smallest, while the 2-ranks
$A_{3\sigma+4}(d)$ for the Glynn type (II) case are the largest. 
This observation is established in Lemma~\ref{lem:2-rank}.

\medskip
Our first goal is to determine $B_6(d)$ explicitly
in terms of Fibonacci numbers. The result,
Theorem~\ref{thm:number-solutions}, was conjectured by Xiang \cite{qx}.

Given any 
integer $x$ not divisible by $2^d-1$, 
recall that by Theorem~\ref{thm5} with $q=2^d$, $s(x)$
equals the number of $1$'s in the binary representation of $x$ reduced
modulo $2^d-1$.
Now, the binary representation of $2x$ mod $2^d-1$ is simply
the rotation of the binary representation of $x$ mod $2^d-1$ by
``one step'', i.e., all
the digits are moved one step to the left, with a 1 that would be moved
beyond the leftmost place moved (``rotated'') to the rightmost place.
This implies that $s(x)=s(2x)$. 
Therefore, given any solution to
\eqref{eq:keyeq1},
multiples of the solution 
by powers of 2 will also give solutions to \eqref{eq:keyeq1}.

\begin{Proposition}\label{prop:solutions}
Let $d$ be an integer $\ge2$.
The binary representations of all the solutions $a$ {\em mod} $2^d-1$
to the equation
\begin{equation}\label{eq:keyeq2}
s(a)+s(5a) = s(6a)+1
\end{equation}
can be constructed by the following procedure:

Step 1. Form all the possible binary strings of length $\le d$ by
concatenating blocks of the form $01$, $0011$, $00111$, subject to the
following restrictions:

\begin{itemize} 
\item [(A)] In a string of length $<d$, the rightmost block must be
$01$, and the block $00111$ must not occur.
\item [(B)] In a string of length $d$, the block $00111$ occurs exactly
once, namely, as the rightmost block. {\rm (}In particular, strings of
length $d$ can only occur for odd $d$.{\rm )}
\end{itemize}

Step 2. Given a string of length $k$ obtained through Step 1,
append $(d-k)$ $0$'s on the left to form a string of length $d$.

Step 3. All the solutions $a$ {\em mod} $2^d-1$ 
to {\em\eqref{eq:keyeq2}}, or rather their
binary representations, are obtained by forming all possible rotations
of the strings that were constructed in Step~2.

\end{Proposition}

Before we move on to a proof of the proposition, we illustrate this
construction by an example.

\begin{Example}\label{ex1}
\em
Let $d=9$. Then all possible strings formed in Step~1
are (the bars indicate separation of blocks)
\begin{gather} \notag
01,\ 01\sep 01,\ 01\sep01\sep01,\ 01\sep01\sep01\sep01,\ 
0011\sep01,\\ 0011\sep01\sep01,\ 01\sep0011\sep01,\ 01\sep01\sep00111,\ 
0011\sep00111.
\end{gather}

Hence, according to Steps~2 and 3, the solutions $a$ mod $2^d-1$ to
\eqref{eq:keyeq2}, in binary notation, are all the possible rotations
of
\begin{gather} \notag
000000001,\ 000000101,\ 000010101,\ 001010101,\ 
000001101,\\ 000110101,\ 001001101,\ 010100111,\ 
001100111.
\end{gather}
Therefore, the solutions are
$$1, 5, 21, 85, 13, 53, 77, 167, 103,$$
and their multiples by powers of 2,
for a total of $9\cdot 9=81$ solutions.
\end{Example}

\noindent{\bf Proof of Proposition~\ref{prop:solutions}}. 
It is straightforward to check that the numbers obtained by
the proposed construction are indeed solutions to
\eqref{eq:keyeq2}. We leave this to the reader.

It remains to be shown that there cannot be other solutions to
\eqref{eq:keyeq2}. We do this by performing a case-by-case analysis
which rules out all other possibilities.

{}From now on, when we talk about some number $x$ mod $2^d-1$, 
we always assume that
$x$ is reduced mod $2^d-1$ to a number between 0 and $2^d-2$.
When we refer to $a$, $5a$, and $6a$, these numbers are to be
viewed mod $2^d-1$.

Before we start our case-by-case analysis, it is helpful to recall
what we would like to achieve. We want to find each number $a$ such
that the number of $1$'s in the binary representation of $6a$ 
is exactly one less than 
the number of $1$'s in the binary representation of $a$ plus
the number of $1$'s in the binary representation of $5a$.

On the other hand, by carrying out the addition $a+5a=6a$ (mod
$2^d-1$), 
we see that the number of $1$'s in the binary representation
of $6a$ can be at most the total of 
the number of $1$'s in the binary representations of $a$ and $5a$ minus
the number of instances of a 1 occurring at the same place in $a$ and $5a$. To
give an example, if $d=5$, and if $a=11010$, so $5a=00110$, then the
number of $1$'s in $6a$ can be at most $3+2-1=4$. (The number of $1$'s in
the binary representation of $a$ is 3, the number of $1$'s in the binary
representation of $5a$ is 2, and there is a 1 in both $a$ and $5a$ in
the second position from the right.) In fact, $6a=00001$, so the
number of $1$'s in $6a$ is even less. 

Thus we can limit our search to $a$'s for which there is {\em exactly one}
such instance of a 1 in the same position in $a$ and $5a$. In
particular, if we detect {\em two} such instances in $a$ and $5a$,
then $a$ cannot be a solution to \eqref{eq:keyeq2}. To
guarantee that the number of $1$'s in the binary representation of $6a$
is {\em exactly} one less (as opposed to two or more less), 
we must have a 0 in both $a$
and $5a$ in the position which is immediately to the left of
the position where 1 occurs in both $a$
and $5a$.

In the analysis that follows, it does not matter
whether we consider $a$ or $2^ma$, for a positive integer $m$. In terms of
binary representations this means that it does not matter whether we
consider the binary representation of $a$ or any cyclic permutation of
it.

Now we turn to our case-by-case analysis.
Since Proposition~\ref{prop:solutions} is easily
checked for the cases $d=2,3,4,5$, 
we may assume in the following that $d\ge6$.

\medskip
{\em Case 1. The block $11111$ occurs in the binary representation of
$a$.}
Suppose that $11111$ occurs in the binary representation of $a$. Then,
since $5a=a+4a$, we have
\begin{equation*} 
\begin{matrix} \hphantom{4}a=&\Punkt\Punkt\Punkt\Punkt11111\Punkt\Punkt\\
4a=&\Punkt\Punkt11111\Punkt\Punkt\Punkt\Punkt\\
\multispan{2}{\leaders\hrule height1pt\hfill}\\
5a=&\Punkt\Punkt\Punkt\Punkt11\Punkt\Punkt\Punkt\Punkt\Punkt
\end{matrix}
\end{equation*}
So, there are two instances of $1$'s in the same position in $a$ and
$5a$. Hence, $a$ cannot be a solution to \eqref{eq:keyeq2}.

\smallskip
{\em Case 2. The block $1111$ occurs in the binary representation of
$a$.}
Since we already ruled out that $11111$ occurs, we may assume that
there is a 0 before, and a 0 after $1111$. We have
\begin{equation*} 
\begin{matrix} \hphantom{4}a=&\Punkt\Punkt\Punkt\Punkt011110\Punkt\Punkt\\
4a=&\Punkt\Punkt011110\Punkt\Punkt\Punkt\Punkt\\
\multispan{2}{\leaders\hrule height1pt\hfill}\\
5a=&\Punkt\Punkt\Punkt\Punkt\Punkt\Punkt\Punkt\underline\Punkt\Punkt\Punkt\Punkt
\Punkt
\end{matrix}
\end{equation*}
Now, if in the underlined position there is no 1 carrying over from
the right, then the addition would give
\begin{equation*} 
\begin{matrix} \hphantom{4}a=&\Punkt\Punkt\Punkt\Punkt011110\Punkt\Punkt\\
4a=&\Punkt\Punkt011110\Punkt\Punkt\Punkt\Punkt\\
\multispan{2}{\leaders\hrule height1pt\hfill}\\
5a=&\Punkt\Punkt\Punkt\Punkt\Punkt101\Punkt\Punkt\Punkt\Punkt
\end{matrix}
\end{equation*}
Again we are encountering two positions with $1$'s in both $a$ and $5a$, and
thus $a$ cannot be a solution to \eqref{eq:keyeq2}.

On the other hand, if a 1 does carry over, then we have
\begin{equation*} 
\begin{matrix} \hphantom{4}a=&\Punkt\Punkt\Punkt\Punkt011110\Punkt\Punkt\\
4a=&\Punkt\Punkt011110\Punkt\Punkt\Punkt\Punkt\\
\multispan{2}{\leaders\hrule height1pt\hfill}\\
5a=&\Punkt\Punkt\Punkt\Punkt\Punkt110\Punkt\Punkt\Punkt\Punkt
\end{matrix}
\end{equation*}
So again, $a$ cannot be a solution to \eqref{eq:keyeq2}.

\smallskip
{\em Case 3. The block $111$ occurs in the binary representation of
$a$.}
Since we already ruled out that $1111$ occurs, we may assume that
there is a 0 before, and a 0 after $111$. We have
\begin{equation*} 
\begin{matrix} \hphantom{4}a=&\Punkt\Punkt\Punkt\Punkt01110\Punkt\Punkt\\
4a=&\Punkt\Punkt01110\Punkt\Punkt\Punkt\Punkt\\
\multispan{2}{\leaders\hrule height1pt\hfill}\\
5a=&\Punkt\Punkt\Punkt\Punkt\Punkt\Punkt\Punkt\Punkt\Punkt\Punkt\Punkt
\end{matrix}
\end{equation*}

If to the left of the block $01110$ we would have a 1,
\begin{equation*} 
\begin{matrix} \hphantom{4}a=&\Punkt\Punkt\Punkt\Punkt101110\Punkt\Punkt\\
4a=&\Punkt\Punkt101110\Punkt\Punkt\Punkt\Punkt\\
\multispan{2}{\leaders\hrule height1pt\hfill}\\
5a=&\Punkt\Punkt\Punkt\Punkt\Punkt\Punkt\Punkt\underline\Punkt\Punkt\Punkt\Punkt
\Punkt
\end{matrix}
\end{equation*}
then, regardless of whether a 1 carries over to the underlined
position or not (see Case~2 for an analogous consideration), we would
obtain two positions with $1$'s in both $a$ and $5a$, and so $a$ cannot 
be a solution to \eqref{eq:keyeq2}.

Therefore, the only way to have the block $111$ is in the form 
$\dots001110\dots$.

\smallskip
{\em Case 4. The block $0110$ occurs in the binary representation of
$a$.}
Let us first assume that we have a 1 on the left of this block, 
so that we encounter
\begin{equation*} 
\begin{matrix} \hphantom{4}a=&\Punkt\Punkt\Punkt\Punkt10110\Punkt\Punkt\\
4a=&\Punkt\Punkt10110\Punkt\Punkt\Punkt\Punkt\\
\multispan{2}{\leaders\hrule height1pt\hfill}\\
5a=&\Punkt\Punkt\Punkt\Punkt\Punkt\Punkt\underline\Punkt\Punkt\Punkt\Punkt\Punkt
\end{matrix}
\end{equation*} 
Now, if no 1 carries over to the underlined position, then we have
\begin{equation*} 
\begin{matrix} \hphantom{4}a=&\Punkt\Punkt\Punkt\Punkt10110\Punkt\Punkt\\
4a=&\Punkt\Punkt10110\Punkt\Punkt\Punkt\Punkt\\
\multispan{2}{\leaders\hrule height1pt\hfill}\\
5a=&\Punkt\Punkt\Punkt\Punkt\Punkt1\underline1\Punkt\Punkt\Punkt\Punkt
\end{matrix}
\end{equation*} 
So, there is a 1 in that position in both $a$ and $5a$, and there is a
1 just left of that position in $5a$. By what we observed earlier, such a
number $a$ cannot be a solution to \eqref{eq:keyeq2}. 
If 1 does carry over to the underlined position, then we have
\begin{equation*} 
\begin{matrix} \hphantom{4}a=&\Punkt\Punkt\Punkt x10110\Punkt\Punkt\\
4a=&\Punkt\Punkt10110\Punkt\Punkt\Punkt\Punkt\\
\multispan{2}{\leaders\hrule height1pt\hfill}\\
5a=&\Punkt\Punkt\Punkt y10\underline0\Punkt\Punkt\Punkt\Punkt
\end{matrix}
\end{equation*} 
where either $x$ or $y$ equals 1. Again, $a$ cannot be a solution to
\eqref{eq:keyeq2}. 

Therefore, the only way to have a block $0110$ is in the form
$\dots00110\dots$. 

\medskip
At this point, let us summarize our observations. The conclusions from
the considerations in Cases~1--4 are that the (cyclic) 
binary strings of length $d$ for $a$
that we are looking for must be built out of the blocks $00111$,
$0011$, and $01$, with possibly filling 0's placed in between.

It is not difficult to see that it is only possible to insert filling
0's in one place in a (cyclic) binary string. For, otherwise we would have at
least two instances of a 1 in the same place in $a$ and $5a$. 
For the same reason, if a block
$00111$ occurs, then it can appear only once, and no filling 0's are
allowed. Finally, it is easy to see that 
an all 0 string, or a string built out of blocks of the form
$0011$ plus filling 0's, or a string
where filling 0's are adjacent to two consecutive $1$'s,
does not give a solution to
\eqref{eq:keyeq2}. All the remaining possibilities are indeed covered
by the construction in the statement of
Proposition~\ref{prop:solutions}. This concludes the proof of the
proposition. 
\QED

Proposition~\ref{prop:solutions} yields the following
formula for the number $B_6(d)$ of solutions to \eqref{eq:keyeq2}.

\begin{Theorem} \label{thm:number-solutions}
Write $A_6(d)=B_6(d)/d$, where 
$B_6(d)$ is the number of solutions mod $2^d-1$ to equation {\em
\eqref{eq:keyeq2}}. Then for $d\ge6$, $A_6(d)$ satisfies the recurrence
\begin{equation} \label{eq:recurrence1}
{A_6(d)} ={A_6(d-2)} +{A_6(d-4)} +1,
\end{equation}
with initial conditions $A_6(2)=0$, $A_6(3)=1$, $A_6(4)=1$, $A_6(5)=3$. 
Equivalently, for any positive integer $m$ we have 
$B_6(2m)=2m(F_m-1)$ and $B_6(2m+1)=(2m+1)(2F_m-1)$, where $F_n$ is the
$n$-th Fibonacci number {\em(}$F_0=F_1=1${\em)}.
\end{Theorem}

{\Proof}
Clearly, the recurrence~\eqref{eq:recurrence1} can be viewed as a
recurrence for the number $B_6(d)/d$ of strings that are constructed in Step~1
of Proposition~\ref{prop:solutions}, i.e., before appending 0's and rotating
the obtained strings of length $d$. 
(There are $B_6(d)/d$ such strings, since due to the particular form of
the strings, all the $d$ 
rotations in Step~3 are indeed different from each other.)

On the basis of this observation, it is easy to demonstrate the
recurrence~\eqref{eq:recurrence1}. For, given $d$, the strings
that are constructed in Step~1 can be separated into three sets:
first, there is the string $01$ of length $2$, 
then there is the set of strings of
length greater than $2$ with leftmost block $01$, and finally there is the
set of strings with leftmost block $0011$. The former set of strings can be
obtained by performing Step~1 with $d$ replaced by $d-2$ and then
appending $01$ to the left of each of the obtained strings. 
The latter set of strings can be
obtained by performing Step~1 with $d$ replaced by $d-4$ and then
appending $0011$ to the left of each of the obtained strings. 
Summing, we obtain the 
recurrence~\eqref{eq:recurrence1} for $B_6(d)/d$. 

{}From the recurrence, the explicit formulas for
$B_6(d)$ in terms of Fibonacci numbers follow immediately. 
\QED

As a corollary, we have the following result on the 2-ranks of those
cyclic difference sets constructed from the Segre hyperovals. 

\begin{Corollary}\label{cor1}
Let $d$ be an odd integer, $d\geq 5$, and let $D_{6,d}$ be the
$(2^d-1,2^{d-1}-1,2^{d-2}-1)$ cyclic difference set in $\Ff_{2^d}^*$
corresponding to the Segre hyperoval $D(x^6)$,
as in Theorem~{\em\ref{thm2}} with $k=6$. Let $\overline{D_{6,d}}$ be the
complement of $D_{6,d}$ in $\Ff_{2^d}^*$. Then the $2$-rank $B_6(d)$
of $\overline{D_{6,d}}$ is equal to $d(2F_{ (d-1)/ {2}}-1)$, where
$F_n$ is the $n$th Fibonacci number {\em(}$F_0=F_1=1${\em)}. 
\end{Corollary} 

{\Proof} This is immediate from Theorem~\ref{thm7} and
Theorem~\ref{thm:number-solutions}.
\QED

Next we turn to the computation of the 
2-ranks of the cyclic difference set arising
from the two types of Glynn hyperovals. According to
Theorem~\ref{thm7}, \eqref{Glynn-I-oval} and \eqref{Glynn-II-oval},
we need to compute the number $B_k(d)$ of solutions $a$ mod $2^d-1$ to  
\begin{equation} \label{eq:keyeq1'}
s(a)+s((k-1)a) = s(ka)+1,
\end{equation}
where $d$ is odd, $d\geq 7$, and $k=\sigma+\gamma$ or $k=3\sigma +4$, with
$\sigma=2^{ (d+1)/ {2}}$, and $\gamma=2^{(3d+1)/4}$ if $d\equiv1$ mod
4, whereas $\gamma=2^{(d+1)/4}$  
if $d\equiv3$ mod 4.
These are much more difficult problems than the counting
problem in the case $k=6$. We resolve them 
in Theorems~\ref{conj3} and \ref{conj2} by making appeal to the
so-called ``transfer matrix method" (see \cite[Sec.~4.7]{Stanley}).
(We could have proved Theorem~\ref{thm:number-solutions} by means of
the transfer matrix method as well. We preferred to take the more
direct approach via Proposition~\ref{prop:solutions}.)

In the proofs of Theorems~\ref{conj3} and \ref{conj2}, we will 
make use of the result (in the folklore) that when it is known (by some
abstract means) that a sequence satisfies {\em some} linear
recurrence, and when a bound for the order of the recurrence is also known,
then one needs only to check a certain number of special instances of a {\em
specific} recurrence to prove that the sequence satisfies this
recurrence ``always". We give the precise statement in the lemma
below. For the sake of completeness, we also supply a proof.

\begin{Lemma} \label{lem:rec-always}
Let
$(f_n)_{n\ge0}$ be a sequence of complex numbers. Suppose that we
know that the generating function $\sum _{n\ge0} ^{}f_nz^n$ for the
sequence is rational, i.e., that it equals $p(z)/q(z)$, where $p(z)$ and $q(z)$
are polynomials in $z$, and that the degree of the numerator, $p(z)$,
is at most $P$, and the degree of the
denominator, $q(z)$, is at most $Q$. If the sequence $(f_n)$
satisfies the recurrence
\begin{equation} \label{eq:rec}
\sum_{i=0}^k a_if_{n-i} = c
\end{equation}
for $n=n_0,\dots,N$, where $n_0 \ge k $, $N=\max\{P+k+1,Q+n_0\}$, and where
$a_0,a_1\dots,a_k$ and $c$ are some given complex numbers, then the
recurrence {\em\eqref{eq:rec}} is satisfied for all $n \ge n_0$.
\end{Lemma}

{\Proof}
Since \eqref{eq:rec} is satisfied for $n=n_0,\dots,N$, we have 
\begin{equation} \label{eq:formal}
 a(z)\sum _{n\ge0} ^{}f_nz^n= r(z) + c/(1-z) + z^{N+1}S(z),
\end{equation}
where $a(z)=\sum_{i=0}^ka_iz^i$, $r(z)$ is a polynomial of degree  
less than $n_0$, and $S(z)$ is some formal power series, the coefficients
of which ``measure" how \eqref{eq:rec} fails for $n>N$. In particular,
$S(z)$ vanishes if and only if  \eqref{eq:rec} holds for {\em all\/} 
$n \ge n_0$.
Now we multiply both sides of \eqref{eq:formal} by $(1-z)q(z)$.
This gives
$$
(1-z)a(z)p(z) = (1-z)q(z)r(z) + cq(z) +
z^{N+1}(1-z)q(z)S(z). 
$$
In this equation, the term on the left-hand side and the first two
terms on the right-hand side all have degree less than $N+1$. It
follows that $z^{N+1}(1-z)q(z)S(z)=0$, whence $S(z)=0$.
This completes the proof.
\QED

Now we are ready to determine the number of solutions to
\eqref{eq:keyeq1'} for the Glynn cases.
We start with the Glynn type (II) case because it is easier to handle.

\begin{Theorem}\label{conj3}
Let $d$ be an odd integer, $d\geq 3$. Let $B_{3\sigma+4}(d)$ be
defined as in the first paragraph of this section, 
and let $A_{3\sigma+4}(d)=B_{3\sigma+4}(d)/d$. Then 
\begin{equation}\label{rec:Glynn2}
A_{3\sigma+4}(d)=A_{3\sigma+4}(d-2)+3A_{3\sigma+4}(d-4)-
A_{3\sigma+4}(d-6)-A_{3\sigma+4}(d-8)+1,
\end{equation}
for all odd $d$, $d\geq 11$,
with initial values as given in Table~{\em1}.
\end{Theorem}

{\Proof}
First of all, from the proof of
Proposition~\ref{prop:solutions}, we see that the solutions $a$, $0<a<2^d-1$,
to equation \eqref{eq:keyeq1'} are characterized by the following
property (the pertinent considerations in the proof of
Proposition~\ref{prop:solutions} are valid for generic $k$, i.e., 
they do not restrict to the case $k=6$):
\begin{enumerate}
\item[({\refstepcounter{equation}\label{no1}}\ref{no1})]
There is exactly one instance of a 1 occurring at the same place in
$a$ and $(k-1)a$, and immediately 
to the left of those 1's there is a 0 in both $a$ and $(k-1)a$,
\end{enumerate}
where we view
$a$ and $(k-1)a$ (mod $2^d-1$) as binary strings of length $d$.

We are actually concerned not with the number $B_k(d)$ of {\em all\/}
solutions to \eqref{eq:keyeq1'}, but rather with $A_k(d)=B_k(d)/d$. Since any 
solution $a$ immediately gives 
rise to $d$ solutions by forming all possible
rotations of $a$, we may, without loss of generality, restrict our
attention to those solutions $a$ where the instance of a 1 imposed by
\eqref{no1} 
occurs exactly in the {\em middle} of the string.
The number $A_k(d)$ then is the number of these {\em special\/}
solutions.

Now, for a proof of the theorem, we define a bijection between
these special solutions to \eqref{eq:keyeq1'} (with $k=3\sigma+4$, where
$\sigma=2^{(d+1)/2}$) and 
a set of closed walks in a certain directed
graph $D$. As we demonstrate at the end, it then follows from first
principles that the generating function for the number of these closed walks, 
and hence for the number of solutions to \eqref{eq:keyeq1'}, is rational (see 
also
\cite[Sec.~4.7]{Stanley}). Finally, we check \eqref{rec:Glynn2} for
enough (odd) values of $d$ so that Lemma~\ref{lem:rec-always} will imply
that \eqref{rec:Glynn2} has to hold for {\em all\/} odd $d\ge11$.

In order to motivate the construction of the aforementioned directed
graph $D$, we analyse, for a given special solution $a$ to \eqref{eq:keyeq1'}, 
the computation of $(k-1)a$ and $ka$ mod $2^d-1$. Let $d=2r-1$, with 
$r \ge 2$, and let the
binary representation of $a$ be $a_{2r-2}\dots a_1a_0$. Since
$k-1=3\sigma+3=3\cdot 2^{(d+1)/2}+3=2^{r+1}+2^r+2+1$, we may compute $(k-1)a$
and $ka=(k-1)a+a$ mod $2^d-1$ as follows:
\begin{equation} \label{eq:add}
\begin{matrix} 
\hfill a=&a_{2r-2}&\dots&a_r&a_{r-1}&a_{r-2}&\dots&a_0\\
\hfill 2a=&a_{2r-3}&\dots&a_{r-1}&a_{r-2}&a_{r-3}&\dots&a_{2r-2}\\
\hfill 2^{r+1}a=&a_{r-3}&\dots&a_{2r-2}&a_{2r-3}&a_{2r-4}&\dots&a_{r-2}\\
\hfill 2^ra=&a_{r-2}&\dots&a_0&a_{2r-2}&a_{2r-3}&\dots&a_{r-1}\\
\omit&\multispan{7}{\leaders\hrule height1pt\hfill}\\
\hfill (k-1)a=&y_{2r-2}&\dots&y_r&y_{r-1}&y_{r-2}&\dots&y_0\\
\hfill a=&a_{2r-2}&\dots&a_r&a_{r-1}&a_{r-2}&\dots&a_0\\
\omit&\multispan{7}{\leaders\hrule height1pt\hfill}\\
\hfill ka=&z_{2r-2}&\dots&z_r&z_{r-1}&z_{r-2}&\dots&z_0\\
\end{matrix}
\kern-7.5cm
\raise-1.1cm\vbox{\hsize1cm
$$\sideset {}+
     {\left.\vbox{\hsize0pt\vskip1cm}\right\}}$$
\vskip-6pt
$$\sideset {}+
     {\left.\vbox{\hsize0pt\vskip0.5cm}\right\}}$$}
\end{equation}
where
\begin{align} \label{eq:add1}
a_i+a_{i-1}+a_{r+i-2}+a_{r+i-1}+b_{i-1}&=y_i+2b_i, \quad 
\text {for }i=0,1,\dots,2r-2,\\
\label{eq:add2}
y_i+a_i+c_{i-1}&=z_i+2c_i,\quad 
\text {for }i=0,1,\dots,2r-2,
\end{align}
for some integers $b_i$, $c_i$ with $0 \le b_i \le 3$, $i=0,1,\dots,2r-2$,
$c_{r-1}=1$, and $c_i=0$ for $i=0,1,\dots,r-2,r,\dots,2r-2$,
with the convention that indices are read modulo $2r-1$ throughout.
The integer $b_i$ is the carry-over from the $i$th
place to the $(i+1)$th 
during the addition that gives $(k-1)a$ in \eqref{eq:add}, while 
$c_i$ is the carry-over from the $i$th
place to the $(i+1)$th 
during the addition of $(k-1)a$ and $a$ in \eqref{eq:add} that gives $ka$.
Here, places are counted from the right starting with 0; for example, 
in the top row of
\eqref{eq:add}, the digit $a_0$ is in the 0th place and
the digit $a_1$ is in the 1th place.
(Indeed, as is not difficult to see, there is a unique such choice
of $b_i$'s and $c_i$'s which satisfy \eqref{eq:add1} and
\eqref{eq:add2}
such that the $y_i$'s and $z_i$'s are the digits in the
additions given by \eqref{eq:add}.)
Equation \eqref{eq:add2}, together with the above
assignments for the $c_i$'s,
reflects the conditions imposed by \eqref{no1} with the 1's in the
middle place.

The subsequent considerations become more transparent, if in 
\eqref{eq:add1} and \eqref{eq:add2} we perform a relabelling. Given a
sequence $(x_i)_{i=0,\dots,2r-2}$, 
where $(x_i)$ can be any of $(a_i)$, $(b_i)$, $(c_i)$, $(y_i)$, or
$(z_i)$, define the new sequence $(\tx_j)_{j=0,\dots,2r-2}$ 
by $\tx_j:=x_{(r-1)j}$.
Note that since $(r-1,d)=(r-1,2r-1)=1$, 
the sequence $(\tx_j)_{j=0,\dots,2r-2}$ is a permutation of
$(x_i)_{i=0,\dots,2r-2}$ 
(recall that indices are read modulo $2r-1$). 
Now, by putting $i=(r-1)j$, equations 
\eqref{eq:add1} and \eqref{eq:add2} read, in terms of the transformed 
sequences, 
\begin{align} \label{eq:tadd1}
\ta_j+\ta_{j+2}+\ta_{j+3}+\ta_{j+1}+\tb_{j+2}&=\ty_j+2\tb_j, \quad 
\text {for }j=0,1,\dots,2r-2,\\
\label{eq:tadd2}
\ty_j+\ta_j+\tc_{j+2}&=\tz_j+2\tc_j,\quad 
\text {for }j=0,1,\dots,2r-2,
\end{align}
where
$\tb_j$, $\tc_j$ are integers with $0 \le \tb_j \le 3$, $j=0,1,\dots,2r-2$,
$\tc_1=1$, and $\tc_j=0$ for $j=0,2,3,\dots,2r-2$, 
where again all indices are viewed modulo $2r-1$.

Construct a directed graph $D$ as follows.
The set of vertices of $D$ is the set of all vectors
$(a',a'',a''',b',b'',c',c'')$ with $0\le
a',a'',a'''\le 1$, $0\le b',b''\le 3$, and $0\le c',c''\le 1$.

Corresponding to any special solution $a$, $0<a<2^d-1$,
to equation \eqref{eq:keyeq1'} with $k=3\sigma+4$,
define (in the notation of \eqref{eq:tadd1}--\eqref{eq:tadd2}) 
the $2r-1$ vertices
\begin{equation} \label{eq:vertex1}
v_j =  (\ta_{j},\ta_{j+1}, \ta_{j+2},\tb_{j},\tb_{j+1},
\tc_{j},\tc_{j+1}),
\end{equation}
$0 \leq j \leq 2r-2$.
Our construction of the directed edges will be motivated by the desire
to obtain, via equations 
\eqref{eq:tadd1} and \eqref{eq:tadd2},
a bijection between the special solutions $a$
and certain closed walks in $D$ of length $d=2r-1$.
A solution $a$ will correspond to the walk
\begin{equation} \label{eq:walk1}
v_0\ \to\ v_1\ \to\ v_2\ \to\cdots\to\ v_{2r-2}\ \to\ v_0,
\end{equation}
and vice versa.

We connect a vertex
$(a',a'',a''',b',b'',c',c'')$
to a vertex $(A',A'',A''',B',B'',C',C'')$ by a directed edge if and only if
\begin{align} \label{eq:edgeA}
a''=A',&\\
\label{eq:edgeA1}
a'''=A'',&\\
\label{eq:edgeA2}
b''=B',&\\
\label{eq:edgeA3}
c''=C',&\\
\label{eq:edgeA4}
Y:=(a'+a'''+A'''+a''+B''-2b')&\in\{0,1\},\\
\label{eq:edgeA5}
(Y+a'+C''-2c')&\in\{0,1\}.
\end{align}
To motivate the construction of such directed edges,
the reader should think of $(a',a'',a''',\break b',b'',c',c'')$ and
$(A',A'',A''',B',B'',C',C'')$ as
$v_j=(\ta_{j},\ta_{j+1}, \ta_{j+2},\tb_{j},\tb_{j+1},\tc_j,\tc_{j+1})$ and
$v_{j+1}=(\ta_{j+1},\ta_{j+2},\ta_{j+3},\tb_{j+1},\tb_{j+2},\tc_{j+1},\tc_{j+2})$, 
respectively, $j=0,1,\dots,2r-2$, and should
observe that \eqref{eq:edgeA4} reflects \eqref{eq:tadd1}, and that
\eqref{eq:edgeA5} reflects \eqref{eq:tadd2}.

Define the sets $V_0,V_1,V_2$ as the sets of vertices
$(a',a'',a''',b',b'',c',c'')$ for which $(c',c'')$ equals $(0,0)$,
$(0,1)$, $(1,0)$, respectively.
With these definitions, we claim that for odd $d\ge3$
the special solutions $a$, $0<a<2^d-1$, 
to equation \eqref{eq:keyeq1'} with $k=3\sigma+4$
are in bijection with the {\em special} closed walks of length $d$ 
which start in $V_1$, in the first step move to $V_2$, and from there on
visit only vertices from $V_0$, until they return to the starting
vertex from $V_1$.
We leave it to the reader to check the claim in detail, as this
is a straightforward task (cf\@. \eqref{eq:walk1} and
\eqref{eq:vertex1}).

Finally, let $\AM_{ij}$ be the adjacency matrix of $D$ restricted 
to the information about edges from $V_i$ to $V_j$, $i,j\in
\{0,1,2\}$, that is,
$\AM_{ij}$ is the matrix with rows labelled by the vertices from $V_i$ and
columns labelled by the vertices from $V_j$, with an entry 1 in some
row and column whenever
that particular vertex of $V_i$ is connected to that particular
vertex of $V_j$, all other entries being 0.
Then, by basic principles (which are easily verified from scratch),
the number of the above defined special closed walks 
of length $d$ is equal
to the trace of $\AM_{12}\AM_{20}\AM_{00}^{d-3}\AM_{01} $. 
Hence, by means of our bijection between
special solutions to \eqref{eq:keyeq1'} and special closed walks in $D$, 
we have
\begin{align} \notag
\sum_{r \geq 2} A_{3\cdot 2^r+4}(2r-1) z^{r-2}&=\tr\bigg( \sum _{j=0}
^{\infty}\AM_{12}\AM_{20}\AM_{00}^{2j}\AM_{01}z^{j}\bigg)\\
\label{eq:rational}
&=\tr \bigg(\AM_{12}\AM_{20}(\IM-\AM_{00}^2 z)^{-1}\AM_{01} \bigg) .
\end{align}
We would be done if we could compute the expression
on the right-hand side of \eqref{eq:rational}. 
However, the matrix $\AM_{00}$ is a $128\times 128$ matrix, and with
today's technology it is still not easy to compute the inverse of a
$128\times 128$ matrix, such as $(\IM-\AM_{00}^2 z)$, which has the indeterminate
$z$ sitting inside. We circumvent this difficulty by appealing to
Lemma~\ref{lem:rec-always}.

On the right-hand side of \eqref{eq:rational} we have the trace of a matrix,
all the entries of which are rational in $z$. Hence, the right-hand
side, and so also the left-hand side, is rational in $z$. 
Taking into account that the matrix $\AM_{00}$ is a $128\times 128$ matrix, 
we see that we may write the left-hand side 
as $p(z)/q(z)$ where the degree of $p(z)$
is at most $127$ and the degree of $q(z)$ is at most $128$. Now we
are in the position to invoke Lemma~\ref{lem:rec-always}, 
with $n_0 = k = 4$, $P = 127$, $Q = 128$,
once we are able to verify the recurrence \eqref{rec:Glynn2} for 
$d=11, 13 ,\dots,267$. In view of our characterization of 
$A_{3\sigma+4}(d)=A_{3\cdot 2^r+4}(2r-1)$ as the trace of 
$\AM_{12}\AM_{20}\AM_{00}^{d-3}\AM_{01}$, this is 
routinely performed on a 
computer\footnote{{\sl Maple} and {\sl Mathematica} codes for 
generating the matrices 
$\AM_{00}$, $\AM_{01} $, $\AM_{12}$, $\AM_{20} $ can be obtained on
request from the authors, or by WWW at {\tt
http://radon.mat.univie.ac.at/People/kratt/artikel/glynn.html}. The
verification of \eqref{rec:Glynn2} for $d=11,13,\dots,267$
took 8~minutes on a Pentium with 133~MHz.}.

Thus, the proof of the theorem is complete.
\QED


By combining Theorems~\ref{conj3} and \ref{thm7},
we are now
able to precisely determine the 2-ranks of the difference sets corresponding
to the Glynn type (II) hyperovals.

\begin{Corollary}\label{cor3}
Let $d$ be an odd integer, $d\geq 3$, and set $\sigma=2^{(d+1)/2}$. 
Let $D_{3\sigma+4,d}$ be the
$(2^d-1,2^{d-1}-1,2^{d-2}-1)$ cyclic difference set in $\Ff_{2^d}^*$
corresponding to the Glynn type (II) hyperoval $D(x^{3\sigma +4})$,
as in Theorem~{\em\ref{thm2}}. Let $\overline{D_{3\sigma+4,d}}$ be the
complement of $D_{3\sigma+4,d}$ in $\Ff_{2^d}^*$. Then the
$2$-rank $B_{3\sigma+4}(d)$ 
of $\overline{D_{3\sigma+4,d}}$ equals $dA_{3\sigma+4}(d)$,
where  $A_{3\sigma+4}(3)=1$,  $A_{3\sigma+4}(5)=1$,
$A_{3\sigma+4}(7)=5$,  $A_{3\sigma+4}(9)=7$, and for $d \geq 11$,
$A_{3\sigma+4}(d)$  satisfies the recurrence {\em \eqref{rec:Glynn2}}.
\end{Corollary} 

Next we turn to the Glynn type (I) case.
\begin{Theorem}\label{conj2}
Let $d$ be an odd integer, $d\geq 3$. Let $B_{\sigma+\gamma}(d)$ be
defined as in the first paragraph of this section, and let
$A_{\sigma+\gamma}(d)=B_{\sigma+\gamma}(d)/d$. Then 
\begin{equation}\label{rec:Glynn1}
A_{\sigma+\gamma}(d)=A_{\sigma+\gamma}(d-2)+
A_{\sigma+\gamma}(d-4)+A_{\sigma+\gamma}(d-6)+A_{\sigma+\gamma}(d-8)-1,
\end{equation}
for all odd $d$, $d\geq 13$,
with initial values as given in Table~{\em1}.
\end{Theorem}

{\Proof}
We proceed in the spirit of the
proof of Theorem~\ref{conj3}.  
Here we regard
$A_{\sigma+\gamma}(d)$ as the number of special solutions $a$, 
$0 < a < 2^d-1$, to \eqref{eq:keyeq1'} 
(where $k=\sigma + \gamma$) with the instance of a 1
imposed by \eqref{no1} occurring exactly at the $m$th place if
$d=4m-1$, and exactly at the $(3m+1)$th place if
$d=4m+1$ (where, again,
we count places from the right starting with 0).
We will define a bijection between these special solutions and
certain closed walks in a directed graph $D$.

We distinguish between the cases $d=4m-1$ and $d=4m+1$, noting that
$\gamma$ is defined differently in these two cases (see
\eqref{Glynn-I-oval}). 

First let $d=4m-1$, $m\ge2$.
In this case we
have $\sigma=2^{2m}$ and $\gamma=2^m$, by \eqref{Glynn-I-oval}.

We analyse, for a given special solution $a$ to \eqref{eq:keyeq1'}, 
the computation of $(k-1)a$ and $ka$ mod $2^d-1=2^{4m-1}-1$. Let the
binary representation of $a$ be $a_{4m-2}\dots a_1a_0$. Since
$k=\sigma+\gamma=2^{2m}+2^m$, we may compute $ka$ in the following
two ways:
\begin{equation} \label{eq:Add-1}
\setcounter{MaxMatrixCols}{13}
\begin{matrix} 
\hfill 
2^{2m}a=&a_{2m-2}&\kern-2pt\dots\kern-2pt&a_{m}&a_{m-1}&\kern-2pt\dots\kern-2pt&
a_{0}&a_{4m-2}&\kern-2pt\dots\kern-2pt&a_{3m-1}&a_{3m-2}&\kern-2pt\dots\kern-2pt
&a_{2m-1}\\
\hfill 
2^{m}a=&a_{3m-2}&\kern-2pt\dots\kern-2pt&a_{2m}&a_{2m-1}&\kern-2pt\dots\kern-2pt&
a_{m}&a_{m-1}&\kern-2pt\dots\kern-2pt&a_{0}&a_{4m-2}&\kern-2pt\dots\kern-2pt&a_{
3m-1}\\
\omit&\multispan{12}{\leaders\hrule height1pt\hfill}\\
\hfill 
ka=&z_{4m-2}&\kern-2pt\dots\kern-2pt&z_{3m}&z_{3m-1}&\kern-2pt\dots\kern-2pt&
z_{2m}&z_{2m-1}&\kern-2pt\dots\kern-2pt&z_{m}&z_{m-1}&\kern-2pt\dots\kern-2pt&z_
{0}\\
\end{matrix}
\kern-7.5cm
\raise-.15cm\vbox{\hsize1cm
$$\sideset {}+
     {\left.\vbox{\hsize0pt\vskip0.5cm}\right\}}$$}
\end{equation}
and
\begin{equation} \label{eq:Add-2}
\begin{matrix} 
\hfill 
(k-1)a=&y_{4m-2}&\kern-2pt\dots\kern-2pt&y_{3m}&y_{3m-1}&\kern-2pt\dots\kern-2pt
&
y_{2m}&y_{2m-1}&\kern-2pt\dots\kern-2pt&y_{m}&y_{m-1}&\kern-2pt\dots\kern-2pt&y_
{0}\\
\hfill 
a=&a_{4m-2}&\kern-2pt\dots\kern-2pt&a_{3m}&a_{3m-1}&\kern-2pt\dots\kern-2pt&
a_{2m}&a_{2m-1}&\kern-2pt\dots\kern-2pt&a_{m}&a_{m-1}&\kern-2pt\dots\kern-2pt&a_
{0}\\
\omit&\multispan{12}{\leaders\hrule height1pt\hfill}\\
\hfill 
ka=&z_{4m-2}&\kern-2pt\dots\kern-2pt&z_{3m}&z_{3m-1}&\kern-2pt\dots\kern-2pt&
z_{2m}&z_{2m-1}&\kern-2pt\dots\kern-2pt&z_{m}&z_{m-1}&\kern-2pt\dots\kern-2pt&z_
{0}\\
\end{matrix}
\kern-7.5cm
\raise-.15cm\vbox{\hsize1cm
$$\sideset {}+
     {\left.\vbox{\hsize0pt\vskip0.5cm}\right\}}$$}
\end{equation}
where
\begin{align} \label{eq:Add1}
a_i+a_{m+i}+b_{2m+i-1}&=z_{2m+i}+2b_{2m+i}, \quad 
\text {for }i=0,1,\dots,4m-2,\\
\label{eq:Add2}
y_{2m+i}+a_{2m+i}+c_{2m+i-1}&=z_{2m+i}+2c_{2m+i},\quad 
\text {for }i=0,1,\dots,4m-2,
\end{align}
for some integers
$b_i$, $c_i$ with $0 \le b_i \le 1$, $i=0,1,\dots,4m-2$,
$c_m=1$, and $c_i=0$ for $i=0,1,\dots,m-1,m+1,\dots,4m-2$,
with the convention that indices are read modulo $4m-1$ throughout.
The integer $b_i$ is the carry-over from the $i$th
place to the $(i+1)$th during the addition displayed in \eqref{eq:Add-1} 
that gives $ka$, while $c_i$ is the carry-over from the $i$th
place to the $(i+1)$th during the addition displayed in \eqref{eq:Add-2} 
that gives $ka$. 
(Again, it is not difficult to see that there is a unique such choice
of $b_i$'s and $c_i$'s which satisfy \eqref{eq:Add1} and
\eqref{eq:Add2}
such that the $y_i$'s and $z_i$'s are the digits in the
additions given by \eqref{eq:Add-1} and \eqref{eq:Add-2}.)
Equation \eqref{eq:Add2}, together with the above
assignments for the $c_i$'s,
reflects the conditions imposed by \eqref{no1} with the 1's in the
$m$th place.

Next, as in
the proof of Theorem~\ref{conj3}, we perform an appropriate relabelling 
of indices. Given a
sequence $(x_i)_{i=0,\dots,4m-2}$, 
where $(x_i)$ can be any of $(a_i)$, $(b_i)$, $(c_i)$, $(y_i)$, or
$(z_i)$, define the new sequence $(\tx_j)_{j=0,\dots,4m-2}$ 
by $\tx_j:=x_{mj}$.
Note that since $(m,d)=(m,4m-1)=1$, 
the sequence $(\tx_j)_{j=0,\dots,4m-2}$ is a permutation of
$(x_i)_{i=0,\dots,4m-2}$
(recall that indices are read modulo $4m-1$). Now, 
by putting $i=mj$, equations 
\eqref{eq:Add1} and \eqref{eq:Add2} read, in terms of the transformed 
sequences,
\begin{align} \label{eq:tAdd1}
\ta_j+\ta_{j+1}+\tb_{j-2}&=\tz_{j+2}+2\tb_{j+2}, \quad 
\text {for }j=0,1,\dots,4m-2,\\
\label{eq:tAdd2}
\ty_{j+2}+\ta_{j+2}+\tc_{j-2}&=\tz_{j+2}+2\tc_{j+2},\quad 
\text {for }j=0,1,\dots,4m-2,
\end{align}
where $0 \le \tb_j \le 1$, $j=0,1,\dots,4m-2$,
$\tc_1=1$, and $\tc_j=0$ for $j=0,2,3,\dots,4m-2$, 
where again all indices are viewed modulo $4m-1$.

We now interrupt the proof for the case $d=4m-1$ and turn for a moment to
the case $d=4m+1$, $m\ge1$.
Then we have $\sigma=2^{2m+1}$ and $\gamma=2^{3m+1}$, by \eqref{Glynn-I-oval}.

Again, we analyse, for a given special solution $a$ to \eqref{eq:keyeq1'},
the computation of $(k-1)a$ and $ka$ mod $2^d-1=2^{4m+1}-1$.
Let the binary representation of $a$ be $a_{4m}\dots a_1a_0$.
Since $k=\sigma+\gamma=2^{2m+1}+2^{3m+1}$, we may compute $ka$ in the following
two ways:
\begin{equation} \label{eq:ADD-1}
\setcounter{MaxMatrixCols}{13}
\begin{matrix} 
\hfill 
2^{3m+1}a=&a_{m-1}&\kern-2pt\dots\kern-2pt&a_{0}&a_{4m}&\kern-2pt\dots\kern-2pt&
a_{3m+1}&a_{3m}&\kern-2pt\dots\kern-2pt&a_{2m+1}&a_{2m}&\kern-2pt\dots\kern-2pt&
a_{m}\\
\hfill 
2^{2m+1}a=&a_{2m-1}&\kern-2pt\dots\kern-2pt&a_{m}&a_{m-1}&\kern-2pt\dots\kern-2pt&
a_{0}&a_{4m}&\kern-2pt\dots\kern-2pt&a_{3m+1}&a_{3m}&\kern-2pt\dots\kern-2pt&a_{
2m}\\
\omit&\multispan{12}{\leaders\hrule height1pt\hfill}\\
\hfill 
ka=&z_{4m}&\kern-2pt\dots\kern-2pt&z_{3m+1}&z_{3m}&\kern-2pt\dots\kern-2pt&
z_{2m+1}&z_{2m}&\kern-2pt\dots\kern-2pt&z_{m+1}&z_{m}&\kern-2pt\dots\kern-2pt&z_
{0}\\
\end{matrix}
\kern-7.5cm
\raise-.15cm\vbox{\hsize1cm
$$\sideset {}+
     {\left.\vbox{\hsize0pt\vskip0.5cm}\right\}}$$}
\end{equation}
and
\begin{equation} \label{eq:ADD-2}
\begin{matrix} 
\hfill 
(k-1)a=&y_{4m}&\kern-2pt\dots\kern-2pt&y_{3m+1}&y_{3m}&\kern-2pt\dots\kern-2pt&
y_{2m+1}&y_{2m}&\kern-2pt\dots\kern-2pt&y_{m+1}&y_{m}&\kern-2pt\dots\kern-2pt&y_
{0}\\
\hfill 
a=&a_{4m}&\kern-2pt\dots\kern-2pt&a_{3m+1}&a_{3m}&\kern-2pt\dots\kern-2pt&
a_{2m+1}&a_{2m}&\kern-2pt\dots\kern-2pt&a_{m+1}&a_{m}&\kern-2pt\dots\kern-2pt&a_
{0}\\
\omit&\multispan{12}{\leaders\hrule height1pt\hfill}\\
\hfill 
ka=&z_{4m}&\kern-2pt\dots\kern-2pt&z_{3m+1}&z_{3m}&\kern-2pt\dots\kern-2pt&
z_{2m+1}&z_{2m}&\kern-2pt\dots\kern-2pt&z_{m+1}&z_{m}&\kern-2pt\dots\kern-2pt&z_
{0}
\end{matrix}
\kern-7.5cm
\raise-.15cm\vbox{\hsize1cm
$$\sideset {}+
     {\left.\vbox{\hsize0pt\vskip0.5cm}\right\}}$$}
\end{equation}
where
\begin{align} \label{eq:ADD1}
a_{i}+a_{i-m}+b_{2m+i}&=z_{2m+i+1}+2b_{2m+i+1}, \quad 
\text {for }i=0,1,\dots,4m,\\
\label{eq:ADD2}
y_{2m+i+1}+a_{2m+i+1}+c_{2m+i}&=z_{2m+i+1}+2c_{2m+i+1},\quad 
\text {for }i=0,1,\dots,4m,
\end{align}
for some integers
$b_i$, $c_i$ with $0 \le b_i \le 1$, $i=0,1,\dots,4m$,
$c_{3m+1}=1$, and $c_i=0$ for $i=0,1,\dots,3m,3m+2,\dots,4m$,
with the convention that indices are read modulo $4m+1$ throughout.
The integer $b_i$ is the carry-over from the $i$th
place to the $(i+1)$th during the addition displayed in \eqref{eq:ADD-1} 
that gives $ka$, while $c_i$ is the carry-over from the $i$th
place to the $(i+1)$th during the addition displayed in \eqref{eq:ADD-2} 
that gives $ka$. 
(As in the similar places before, there is a unique such choice
of $b_i$'s and $c_i$'s which satisfy \eqref{eq:ADD1} and
\eqref{eq:ADD2}
such that the $y_i$'s and $z_i$'s are the digits in the
additions given by \eqref{eq:ADD-1} and \eqref{eq:ADD-2}.)
Equation \eqref{eq:ADD2}, together with the above
assignments for the $c_i$'s,
reflects the conditions imposed by \eqref{no1} with the 1's in
the $(3m+1)$th place.

Next, we perform an appropriate relabelling 
of indices. Here, given a
sequence\break $(x_i)_{i=0,\dots,4m}$, 
where $(x_i)$ can be any of $(a_i)$, $(b_i)$, $(c_i)$, $(y_i)$, or
$(z_i)$, define the new sequence $(\tx_j)_{j=0,\dots,4m}$ 
by $\tx_j:=x_{-mj}$. (Note
the difference from the previous relabelling.)
Note that since $(-m,d)=(m,4m+1)=1$, 
the sequence $(\tx_j)_{j=0,\dots,4m}$ is a permutation of
$(x_i)_{i=0,\dots,4m}$. Now,
by putting $i=-mj$, equations 
\eqref{eq:ADD1} and \eqref{eq:ADD2} read, in terms of the transformed 
sequences,
\begin{align} \label{eq:tADD1}
\ta_{j}+\ta_{j+1}+\tb_{j-2}&=\tz_{j+2}+2\tb_{j+2}, \quad 
\text {for }j=0,1,\dots,4m,\\
\label{eq:tADD2}
\ty_{j+2}+\ta_{j+2}+\tc_{j-2}&=\tz_{j+2}+2\tc_{j+2},\quad 
\text {for }j=0,1,\dots,4m,
\end{align}
where $0 \le \tb_j \le 1$, $j=0,1,\dots,4m$,
$\tc_1=1$, and $\tc_j=0$ for $j=0, 2,3,\dots,4m$, 
where again all indices are viewed modulo $4m+1$.

Very conveniently, equations \eqref{eq:tADD1}-\eqref{eq:tADD2} are 
identical with equations
\eqref{eq:tAdd1}-\eqref{eq:tAdd2}, so that from now on we can treat
the cases $d=4m-1$ and $d=4m+1$ simultaneously.

Construct a directed graph $D$ as follows.
The set of vertices of $D$ is the set of all vectors
$(a',a'',b',b'',b''',b'''',c',c'',c''',c'''')$ 
with $0\le a',a'',b',b'',b''',b'''',c',c'',c''',c''''\le 1$.

Corresponding to a special solution $a$, $0 < a < 2^d -1$,
to \eqref{eq:keyeq1'} with $k = \sigma + \gamma$,
define (in the notation of 
\eqref{eq:tAdd1}--\eqref{eq:tAdd2}/\eqref{eq:tADD1}--\eqref{eq:tADD2}) 
the $d$ vertices
\begin{equation} \label{eq:vertex2}
v_j = (\ta_j,\ta_{j+1},\tb_{j-2},\tb_{j-1},\tb_{j},
\tb_{j+1},\tc_{j-2},\tc_{j-1},\tc_{j},\tc_{j+1}),
\end{equation} 
$0 \leq j \leq {d-1}$.
Our construction of the directed edges will be motivated by the desire
to obtain, via equations \eqref{eq:tAdd1} and \eqref{eq:tAdd2},
a bijection between the special solutions $a$ and certain closed walks in
$D$ of length $d$.
A solution $a$ will correspond to the walk
\begin{equation} \label{eq:walk2}
v_0\ \to\ v_1\ \to\cdots\to\ v_{d-1}\ \to\ v_0,
\end{equation}
and vice versa.

We connect a vertex
$(a',a'',b',b'',b''',b'''',c',c'',c''',c'''')$
to a vertex $(A',A'',B',B'',B''',\break B'''',C',C'',C''',C'''')$ 
by a directed edge if and only if
\begin{align} 
\label{eq:EDGEA1}
a''=A',&\\
\label{eq:EDGEA2}
b''=B',&\\
\label{eq:EDGEA3}
b'''=B'',&\\
\label{eq:EDGEA4}
b''''=B''',&\\
\label{eq:EDGEA5}
c''=C',&\\
\label{eq:EDGEA6}
c'''=C'',&\\
\label{eq:EDGEA7}
c''''=C''',&\\
\label{eq:EDGEA8}
Z:=(a'+a''+b'-2B'''')&\in\{0,1\},\\
\label{eq:EDGEA9}
(Z+2C''''-A''-c')&\in\{0,1\}.
\end{align}
The reader should think of
$(a',a'',b',b'',b''',b'''',c',c'',c''',c'''')$
and $(A',A'',B',B'',B''',B'''',\break C',C'',C''',C'''')$ as
$v_j=(\ta_j,\ta_{j+1},\tb_{j-2},\tb_{j-1},\tb_{j},\tb_{j+1},
\tc_{j-2},\tc_{j-1},\tc_{j},\tc_{j+1})$
and
$v_{j+1}=(\ta_{j+1},\break \ta_{j+2},\tb_{j-1},\tb_{j},\tb_{j+1},\tb_{j+2},
\tc_{j-1},\tc_{j},\tc_{j+1},\tc_{j+2})$, respectively,
$j=0,1,\dots,d-1$. Equation \eqref{eq:EDGEA8}
reflects \eqref{eq:tAdd1}/\eqref{eq:tADD1}, and equation
\eqref{eq:EDGEA9} reflects \eqref{eq:tAdd2}/\eqref{eq:tADD2}. 

The remainder of the proof proceeds just like the proof of
Theorem~\ref{conj3}. We define sets $V_0,V_1,V_2,V_3,V_4$ as the sets
of vertices $(a',a'',b',b'',b''',b'''',c',c'',c''',c'''')$ for which
$(c',c'',c''',c'''')$ equals $(0,0,0,0)$, $(0,0,0,1)$, $(0,0,1,0)$,
$(0,1,0,0)$, $(1,0,0,0)$, respectively. As is easy to see, 
for odd $d\ge3$
the special solutions $a$, $0<a<2^d-1$, 
to equation \eqref{eq:keyeq1'} with $k=\sigma+\gamma$
are in bijection with the closed walks of length $d$ 
which start in $V_1$, in the first step move to $V_2$, in the next
step move to $V_3$, then move to $V_4$, and
from there on visit only vertices from $V_0$, until 
they return to the starting vertex from $V_1$.
In an analogous manner we define matrices $\AM_{12}$,
$\AM_{23}$, $\AM_{34}$, $\AM_{40}$,
$\AM_{00}$ and $\AM_{01}$,
and obtain a rational expression for the generating function 
for the sequence $(A_{\sigma+\gamma}(2r+1))_{r\ge2}$ in the form
\begin{align*} 
\sum_{r \geq 2} A_{\sigma+\gamma}(2r+1) z^{r-2}&=\tr\bigg( \sum _{j=0}
^{\infty}\AM_{12}\AM_{23}\AM_{34}\AM_{40}\AM_{00}^{2j}\AM_{01}z^{j}\bigg)\\
\label{eq:rational}
&=\tr \bigg(\AM_{12}\AM_{23}\AM_{34}\AM_{40}
(\IM-\AM_{00}^2 z)^{-1}\AM_{01} \bigg) .
\end{align*}
This time the
matrix $\AM_{00}$, the adjacency matrix of $D$ restricted to $V_0$,
is a $64\times 64$ matrix. Therefore, in order to
invoke Lemma~\ref{lem:rec-always}, it suffices to verify
the recurrence \eqref{rec:Glynn1} for 
$d=13,15,\dots,141$. 
This is again routinely performed on a 
computer\footnote{{\sl Maple} and {\sl Mathematica} codes for 
generating the matrices 
$\AM_{00}$, $\AM_{01}$, $\AM_{12}$,
$\AM_{23}$, $\AM_{34}$, $\AM_{40}$ can be obtained on
request from the authors, or by WWW at {\tt
http://radon.mat.univie.ac.at/People/kratt/artikel/glynn.html}. The
verification of \eqref{rec:Glynn1} for $d=13,15,\dots,141$
took 15~seconds on a Pentium with 133~MHz.}. This establishes
\eqref{rec:Glynn1} for {\em all\/} odd $d$, $d\ge13$.

\medskip
Now the proof of the theorem is complete.
\QED

By combining Theorems~\ref{conj2} and \ref{thm7},
we are now
able to precisely determine the 2-ranks of the difference sets corresponding
to the Glynn type (I) hyperovals.

\begin{Corollary}\label{cor2}
Let $d$ be an odd integer, $d\geq 3$, and set $\sigma=2^{(d+1)/2}$, 
$\gamma=2^{(3d+1)/4}$ if $d\equiv1$ {\em mod 4}, and $\gamma=2^{(d+1)/4}$ 
if $d\equiv3$ {\em mod 4}. Let $D_{\sigma+\gamma,d}$ be the
$(2^d-1,2^{d-1}-1,2^{d-2}-1)$ cyclic difference set in $\Ff_{2^d}^*$
corresponding to the Glynn type (I) hyperoval $D(x^{\sigma +\gamma})$,
as in Theorem~{\em\ref{thm2}}. Let $\overline{D_{\sigma+\gamma,d}}$ be the
complement of $D_{\sigma+\gamma,d}$ in $\Ff_{2^d}^*$. Then the
$2$-rank $B_{\sigma+\gamma}(d)$ 
of $\overline{D_{\sigma+\gamma,d}}$ equals $dA_{\sigma+\gamma}(d)$,
where  $A_{\sigma+\gamma}(3)=1$,  $A_{\sigma+\gamma}(5)=1$,
$A_{\sigma+\gamma}(7)=3$,  $A_{\sigma+\gamma}(9)=7$,
$A_{\sigma+\gamma}(11)=13$, and for $d \geq 13$,  $A_{\sigma+\gamma}(d)$
satisfies the recurrence {\em \eqref{rec:Glynn1}}. 
\end{Corollary}

\section{Binary cyclic codes of difference sets
arising from Segre hyperovals}

The next two propositions give some further applications of 
Proposition~\ref{prop:solutions}.
Let ${\mathcal C}_2(D_{6,d})$ be the binary cyclic code of $D_{6,d}$,
as defined near the beginning of Section~3.
Define $\theta(x)$ (viewed as a polynomial over $\Ff_{2}$) by
$$\theta(x) = \sum_{j=0}^{2^d-2} b_j x^j,$$
where $b_0, b_1, \ldots, b_{2^d -2}$ denotes the characteristic sequence
of $D_{6,d}$. 
By suitably ordering the elements of $D_{6,d}$, we may choose a (fixed)
generator $\alpha$ of  $\Ff_{2^d}^*$ such that
$b_j = 1$ if and only if $\alpha^j \in D_{6,d}$.
The significance of $\theta(x)$ is that its greatest common divisor with
$x^{2^d-1} -1$ is the generator polynomial of the
binary cyclic code  ${\mathcal C}_2(D_{6,d})$.  The zeros of $\theta(x)$
in  $\Ff_{2^d}^*$ (i.e., the zeros of the generator polynomial)
are called the {\em{zeros}} of ${\mathcal C}_2(D_{6,d})$.  The 
remaining elements of  $\Ff_{2^d}^*$ are called the {\em{nonzeros}}
of  ${\mathcal C}_2(D_{6,d})$.  Observe that for such a nonzero $u$,
we have $\theta(u)=1$, since  $D_{6,d}$ is invariant under the $d$
automorphisms of  $\Ff_{2^d}$.

\begin{Proposition}\label{cor4}
Let $d$ be an odd integer, $d\geq 3$, and let $D_{6,d}$ be the
$(2^d-1,2^{d-1}-1,2^{d-2}-1)$ cyclic difference set in $\Ff_{2^d}^*$
corresponding to the Segre hyperoval $D(x^6)$,
as in Theorem~{\em\ref{thm2}} with $k=6$. 
Let ${\mathcal C}_2(D_{6,d})$ be the binary cyclic code of $D_{6,d}$. 
Then the nonzeros of ${\mathcal C}_2(D_{6,d})$ are $1$, $\alpha^{-5a}$, 
where $\alpha$ is the primitive element of $\Ff_{2^d}^*$ defined above, and
$a$ runs through all the solutions of\/ {\em{\eqref{eq:keyeq2}}} given 
explicitly 
in Proposition~{\em{\ref{prop:solutions}}}.
\end{Proposition}

{\Proof} Clearly $1$ is a nonzero of ${\mathcal C}_2(D_{6,d})$ because the 
cardinality of $D_{6,d}$ is odd. For $\mathfrak{p}$ a prime ideal in 
$\Zz[\xi_{2^d-1}]$ lying over
$2$, let $\omega=\omega_{\mathfrak p}$ be the Teichm\"uller character on
$\Ff_{2^d}$.  By definition of $\theta(x)$ and $\alpha$,
we see that $\alpha^i$, $0< i\leq 2^d-2$ is a nonzero of 
${\mathcal C}_2(D_{6,d})$ if and only if
$\omega^{i}\left(D_{6,d}\right)$ (mod $\mathfrak{p}$)
$\neq 0$. From the proof of Theorem~\ref{thm7}, we see that 
for $0 < a \leq 2^d-2$,
$\omega^{-5a}\left(D_{6,d}\right)$ (mod $\mathfrak{p}$)
$\neq 0$ if and only if $s(a)+s(5a)=s(6a)+1$. Therefore the nonzeros of 
${\mathcal C}_2(D_{6,d})$ are $1, \alpha^{-5a}$,  where $a$ runs through the 
solutions of \eqref{eq:keyeq2}.
This completes the proof.
\QED

\begin{Example}
\em

Let $D_{6,9}$ be the
$(2^9-1,2^{8}-1,2^{7}-1)$ cyclic difference set in $\Ff_{2^9}^*$
corresponding to the Segre hyperoval $D(x^6)$,
as in Theorem~\ref{thm2} with $k=6$. Let ${\mathcal C}_2(D_{6,9})$ be the 
binary cyclic code of $D_{6,9}$. By Proposition~\ref{cor4} and
Example~\ref{ex1}, we know 
that the nonzeros of ${\mathcal C}_2(D_{6,9})$ are $1$, $\alpha^{-5\cdot 2^i}$, 
$\alpha^{-25\cdot 2^i}$, $\alpha^{-77\cdot 2^i}$, $\alpha^{-117\cdot 2^i}$, 
$\alpha^{-9\cdot 2^i}$, $\alpha^{-19\cdot 2^i}$, $\alpha^{-7\cdot 2^i}$, 
$\alpha^{-37\cdot 2^i}$, and $\alpha^{-2^i}$, where
$i=0,1,2,\ldots ,8$. A simple computation shows that $\alpha,
\alpha^2, \alpha^3, \ldots , \alpha^{42}$ are zeros of ${\mathcal
C}_2(D_{6,9})$. 
By the BCH-bound \cite[p.~201]{MacWiSl}, the code 
${\mathcal C}_2(D_{6,9})$ is a 21-error correcting code of length 511
and dimension 82 over $\Ff_{2}$, (with minimum distance at least 43).
\end{Example}

Let $d$ be an odd integer, $d\geq3$. Let $a$ be an integer
with $1 \leq a \leq 2^d -2$,
and let $r$ be the smallest positive integer such
that $2^{r+1}a\equiv a$ (mod $2^d-1$). The {\em cyclotomic
coset\/} containing $a$ consists of $\{a, 2a, 2^2a, \ldots
,2^ra\}$. The solutions of \eqref{eq:keyeq2} can be partitioned into
$A_6(d)$ disjoint cyclotomic cosets, each of cardinality $d$. Let $J$
be a complete set of cyclotomic coset representatives of the
solutions of \eqref{eq:keyeq2}. We can give a criterion involving trace to
decide when the equation $x^6+x+\beta=0$, $\beta\in \Ff_{2^d}^*$,
has two distinct solutions (as opposed to no solutions) in $\Ff_{2^d}$.  

\begin{Proposition}\label{cor5}
Let $d$ be an odd integer, $d\geq 3$. Then $x^6+x+\beta=0$,
$\beta\in \Ff_{2^d}^*$, has two distinct solutions in $\Ff_{2^d}$
if and only if $\Tr (\sum_{a\in J}\beta^{5a})=0$, where
$\Tr $ is the trace from $\Ff_{2^d}$ to $\Ff_2$, and $J$ is
defined above. 
\end{Proposition}

{\Proof} By Fourier inversion and the definition of $\theta(x)$,
$$(2^d-1)b_j = \sum_{i=0}^{2^d-2} \theta(\alpha^{-i}) \alpha^{ij},$$
for each $j$ with $0 \leq j \leq2^d-2$.
By Proposition~\ref{cor4}, 
$\theta(\alpha^{-i})$ equals 1 if $i=0$ or $i=5a$, and otherwise
$\theta(\alpha^{-i})$ equals 0, where $a$ runs through the solutions
of  \eqref{eq:keyeq2}.  The result now follows with 
$\beta = \alpha^j$.
\QED

\begin{Example}
\em
Let $d=9$. By Example~\ref{ex1}, $J=\{1,5,21,85,13,53,77,167,103\}$ is
a complete set of cyclotomic coset representatives of the
solutions of \eqref{eq:keyeq2}. Therefore $x^6+x+\beta=0$, $\beta\in
\Ff_{2^9}^*$, has two distinct solutions in $\Ff_{2^9}$ if and
only if $\Tr(\beta+\beta^5+\beta^7+\beta^9+\beta^{19}+\beta^{25}+\beta^{37}+
\beta^{77}+\beta^{117})=0$.

\end{Example}

It is interesting to compare Proposition~\ref{cor5} with the result
(following from Hilbert's Theorem 90) that
$x^2+x+\beta=0$, $\beta\in \Ff_{2^d}^*$
has two solutions in $\Ff_{2^d}$ if and only if $\Tr (\beta)=0$.  

\section{Inequivalence of difference sets arising from hyperovals, 
quadratic residue difference sets, and GMW difference sets}

	In the next theorem, we answer a question raised by Maschietti
\cite{ma}
by showing that the difference sets $D_{k,d}$ arising from the hyperovals in
\eqref{reg-oval}, \eqref{Segre-oval}, \eqref{Glynn-I-oval},
\eqref{Glynn-II-oval} are all inequivalent to each other. 
The difference sets arising from the hyperovals in 
\eqref{reg-oval}--\eqref{trans-oval} are of
course equivalent for the same $d$, since as we showed at the end of
Section~2, they are all Singer difference sets.

We precede the statement of the theorem by the following auxiliary result which
is of independent interest. It relates the magnitudes of the
2-ranks of the difference sets arising from Segre and Glynn
hyperovals.

\begin{Lemma} \label{lem:2-rank}
With $A_k(d)$ as defined at the beginning of Section~{\em4}, we have 
for each odd $d \geq 15$ that
\begin{equation} \label{eq:inequ} 
A_{3\sigma+4} (d)  >  A_{\sigma+\gamma} (d)    >   A_6(d),
\end{equation}
i.e., the $2$-ranks of the difference sets that arise from Glynn
type (II) hyperovals dominate the $2$-ranks of the difference sets that 
arise from Glynn type (I) hyperovals, which in turn dominate the $2$-ranks 
of the difference sets arising from Segre hyperovals.
\end{Lemma}

{\Proof}
Since 
$A_6(d)=2F_{(d-1)/2} - 1$ by Theorem~\ref{thm:number-solutions},
it follows that
$$
A_6(d)= \tfrac {2} {\sqrt5}\left(\tfrac {1+\sqrt5}2\right)^{(d+1)/2} -
\tfrac {2} {\sqrt5}\left(\tfrac {1-\sqrt5}2\right)^{(d+1)/2} - 1.
$$
Note that $(1+\sqrt5)/2=1.6181\dots$

Solving the recurrence \eqref{rec:Glynn1}, we obtain
$$A_{\sigma+\gamma}(d)=
a_1\om_1^{(d+1)/2}+a_2\om_2^{(d+1)/2}+
a_3\om_3^{(d+1)/2}+a_4\om_4^{(d+1)/2}+a_5,$$
where $\om_1,\om_2,\om_3,\om_4$ are 
four distinct complex zeros of the
polynomial $z^4-z^3-z^2-z-1$, and where $a_1,a_2,a_3,a_4,a_5$ are
explicitly known complex numbers. The zero with greatest absolute
value is $\om_1$, say, where $\vert\om_1\vert=1.927561975\dots$

Solving the recurrence \eqref{rec:Glynn2}, we obtain
$$A_{3\sigma+4}(d)=b_1\rho_1^{(d+1)/2}+b_2\rho_2^{(d+1)/2}+
b_3\rho_3^{(d+1)/2}+b_4\rho_4^{(d+1)/2}+b_5,$$
where $\rho_1,\rho_2,\rho_3,\rho_4$ are 
four distinct complex zeros of the
polynomial $z^4-z^3-3z^2+z+1$, and where $b_1,b_2,b_3,b_4,b_5$ are
explicitly known complex numbers. The zero with greatest absolute
value is $\rho_1$, say, where $\vert\rho_1\vert=2.095293985\dots$

{}From these formulas, we see that for odd $d>25$,
$$
A_{3\sigma+4} (d) > \tfrac {1} {3}2^{(d+1)/2}  > A_{\sigma+\gamma} (d)
	>  (1.6181)^{(d+1)/2}   >  A_6(d).          
$$
Appealing to Table~1, we see that \eqref{eq:inequ} holds for odd 
$d \geq 15$.

This completes the proof.
\QED

\begin{Theorem} \label{thm8}
The difference sets $D_{k,d}$ arising from the hyperovals in
{\em\eqref{reg-oval}, \eqref{Segre-oval}, \eqref{Glynn-I-oval},
\eqref{Glynn-II-oval}} are all inequivalent.
\end{Theorem}

{\Proof}  It suffices to show that the 2-ranks of the complements
$\overline{D_{k,d}}$
in question are distinct.  
Since the difference set $D_{k,d}$ arising from a regular hyperoval
is a Singer difference set, the 2-rank of its complement equals $d$, by
Theorem~\ref{thm6}
and the remarks following Theorem~\ref{thm6}.  By
Corollary~\ref{cor1}, the 2-rank of
$\overline{D_{6,d}}$ is $dA_6(d) = d(2F_{(d-1)/2} - 1)$, 
which exceeds $d$ for $d \ge 5$.  This 
show the inequivalence for the cases \eqref{reg-oval} and
\eqref{Segre-oval}.
By Table~1
and Lemma~\ref{lem:2-rank}, the 2-ranks in all the cases 
must be distinct.

This completes the proof.
\QED

Let $Q_d$ denote the quadratic residue cyclic $(2^d-1,2^{d-1}-1,2^{d-2}-1)$
difference set in $\Ff^*_{2^d}$ \cite[p.~244]{ju}.  Thus, 
$2^d -1$ must be a prime, and  $Q_d =
\{\alpha^r:r$ is a square mod $2^d-1$,
$0 < r < 2^d-1\}$, where
$\alpha$ generates the group 
$\Ff^*_{2^d}$.
It is easily checked that $Q_2$  and $Q_3$ are equivalent to 
difference sets arising from the regular hyperoval, while $Q_5$ is equivalent
to the difference set arising from the Segre hyperoval.
We now show that for (odd) $d \ge 7$, the sets $Q_d$ are
inequivalent to difference sets arising from hyperovals.

\begin{Theorem} \label{thm9}
Let $d \ge 7$ be odd. 
Then the quadratic residue difference sets $Q_d$ are in-
equivalent to the difference sets $D_{k,d}$ arising from the
hyperovals in {\em\eqref{reg-oval}--\eqref{Glynn-II-oval}}. 
\end{Theorem}

{\Proof} As $\chi$ runs through the (multiplicative)
nontrivial characters on $\Ff^*_{2^d}$, 
$\chi(Q_d) =  \sum _{\beta\in Q_d} ^{}  \chi(\beta)$
assumes exactly two distinct
values $z$ and $\overline z$, where $z$ is a 
complex number and $\overline z$ is its complex 
conjugate.  If $D$ is any difference set in $\Ff^*_{2^d}$ equivalent to $Q_d$,
then as $\chi$ runs through the nontrivial characters on $\Ff^*_{2^d}$,  
$\chi(D)$ assumes only the values  $z$ and $\overline z$ 
up to a factor of a  $(2^d-1)$-th
root of unity.  Now assume for the purpose of contradiction that the 
aforementioned $D$ is one of the difference sets $D_{k,d}$ in 
\eqref{reg-oval}--\eqref{Glynn-II-oval}.
In view of \eqref{eq2.1}, it suffices to show that as $\chi$ runs through the
nontrivial characters on $\Ff^*_{2^d}$, the principal ideals in
$\Zz[\xi_{2^d-1}]$
generated by the Jacobi sums $J(\chi,\chi^{k-1})$ assume {\em more} than two
distinct values.  Thus, by \eqref{eq3.3a},
it suffices to show that as  $a$   runs through
the values $1,2,\dots,2^d-2$, the expression $c(a):=s(a)+s(a(k-1)) - s(ak)$
assumes at least three distinct values. One checks this as follows.

For $k$ as in \eqref{Glynn-II-oval}, $c(1)$, $c(-1)$, and $c(7)$ are distinct.
For $k$ as in \eqref{Glynn-I-oval}, $c(1)$, $c(-1)$, and $c(\sigma+1)$
are distinct. 
For $k$ as in \eqref{Segre-oval}, $c(1)$, $c(-1)$, and $c(3)$ are distinct.
Finally, for $k$ as in \eqref{reg-oval} or \eqref{trans-oval}, $c(1)$,
$c(-1)$, and $c(k+1)$ are distinct.
This completes the proof.
\QED

Next we consider GMW difference sets \cite{go}.
Suppose that $d$ has the factorization $d=uv$ for integers
$u,v > 1$, and write $q=2^u$ (so that $2^d = q^v$).  Let $r$ be an
integer strictly between 1 and $q-1$ which is relatively prime to
$q-1$ and which is not equal to a power of 2.  Define the GMW
set $$G_{d,v,r} = \{y \in  \Ff_{2^{d}}^* \mid \Tr_{q/2}
((\Tr_{q^v/q}(y))^r) = 0 \},$$
where $\Tr_{q^v/q}$ denotes the trace map from $ \Ff_{q^v}$ to
$\Ff_q$. 
(We did not allow $r$ to be a power of 2, for otherwise
$G_{d,v,r}$ would of course be equivalent to a Singer difference set.)
Let $w$ denote the number of $1$'s in the binary expansion of $r$
(so that $1 < w < u$).  It is known \cite{schwel} that $G_{d,v,r}$
is a  $(2^d -1,2^{d-1}-1,2^{d-2}-1)$ difference set whose complement
has 2-rank $uv^w$.  (For the $p$-rank of general
GMW difference sets in fields of characteristic $p$, see \cite{ab}.)
Theorem~\ref{thm10} below shows that the GMW sets $G_{d,v,r}$
are inequivalent to the difference sets $D_{k,d}$ in
\eqref{reg-oval}--\eqref{Glynn-II-oval}. We will first need four
lemmas.

\begin{Lemma} \label{GMWlem1}
Write $d=uv$ and $q=2^u$ as above.
For each nontrivial multiplicative character $\chi$ on $\Ff_{2^d}^*$
we have
$$\chi(G_{d,v,r}) = \sum _{y\in G_{d,v,r}} ^{}\chi(y) = \frac {1} {2} g(\chi)
g_1({\chi_1}^{ 1/r}) /g_1(\chi_1),$$
where $1/r$ denotes the inverse of $r$ modulo $q-1$,
$\chi_1$ denotes the restriction of $\chi$ to $\Ff_q^*$, and
where $g$ and $g_1$ denote the Gauss sums over  $\Ff_{2^d}^*$
and  $\Ff_q^*$, respectively.
\end{Lemma}

{\Proof}
We have
$$2\chi(G_{d,v,r})=\sum_{y \in  \Ff_{2^d}^*} \chi(y)
(-1)^{\Tr_{q/2} ((\Tr_{q^v/q}(y))^r)}.$$
Write $$T= (q^v -1)/(q-1),$$ and let $\alpha$ be a fixed generator of
$\Ff_{2^d}^*,$ so that every $y \in  \Ff_{2^d}^*$ has the form
$y=\alpha^{iT+j},$
with $0 \leq i < q-1$ , $0 \leq j < T$.  Therefore,
since $\alpha^{T} \in \Ff_q^*$,
$$2\chi(G_{d,v,r})=\sum_{j=0}^{T-1}\chi(\alpha^j)\sum_{i=0}^{q-2}
\chi^{1/r}(\alpha^{irT})(-1)^ {\Tr_{q/2} (\alpha^{irT}
(\Tr_{q^v/q}(\alpha^j))^r)}.$$
Thus $$2\chi(G_{d,v,r})=A+B,$$
where
$$A= \sum_{j=0}^{T-1} \chi(\alpha^j) \overline{\chi}(\Tr_{q^v/q}(\alpha^j))
g_1(\chi_1^{1/r})$$ and
$$B= \sum\limits_{\stackrel{j=0}{\Tr_{q^v/q}(\alpha^j)=0}}^{T-1}
{} \sum_{i=0}^{q-2} \chi(\alpha^{iT+j})=
\sum\limits_{\stackrel{y \in \Ff_{2^d}^*} {\Tr_{q^v/q}(y)=0}} \chi(y).$$

To evaluate $A$, observe that
$$A=g_1(\chi_1^{1/r}) \sum\limits_{\stackrel{y \in \Ff_{2^d}^*}
{\Tr_{q^v/q}(y)=1}} \chi(y).$$
This sum on $y$ is an Eisenstein sum, which is known 
\cite[pp.~391, 400]{be1}
to equal $g(\chi)/g_1(\chi_1)$ or $-g(\chi)/q,$ according as
$\chi_1$ is nontrivial or trivial.  Thus
$$A=\begin{cases} g_1(\chi_1^{1/r})g(\chi)/g_1(\chi_1),&\text{if }
\chi_1 \text{ is nontrivial,}\\
g(\chi)/q,
&\text{if }\chi_1 \text{ is trivial.}\end{cases}$$

To evaluate $B$, observe that $B$ is
an Eisenstein sum which is known
\cite[pp.~389, 391, 400]{be1} to equal
$$B=\begin{cases} 0,
&\text{if }\chi_1 \text{ is nontrivial,}\\
g(\chi)(q-1)/q,
&\text{if }\chi_1 \text{ is trivial.}\end{cases}$$
Adding our evaluations of $A$ and $B$, we obtain the desired result.
\QED

\begin{Lemma} \label{GMWlem2}
Let $v=3$, $u \geq 7$ with $u$ odd, $q=2^u$, and $d=uv$.
For an integer $a$ not divisible by $q^v -1$, let $s(a)$ denote
the number of\/ $1$'s in the binary expansion of the reduction of $a$ modulo
$q^v -1$.
As in {\em{\eqref{Glynn-I-oval}}}, define  $\sigma=2^{(d+1)/2}$ and
$\gamma=2^{(3d+1)/4}$ if $d\equiv1$ {\em mod 4}, whereas $\gamma=2^{(d+1)/4}$
if $d\equiv3$ {\em mod 4}. Then for $a=(q-1)(\gamma+(-1)^{(d-1)/2}),$
$$s(a)+s(a(\sigma +\gamma -1)) - s(a(\sigma+\gamma)) < u.$$
\end{Lemma}

{\Proof}
First suppose that $d\equiv3$ mod 4.  Then it may be checked directly
that $s(a)=u$, $s(a(\sigma +\gamma -1)) = (d-1)/2$, and
$s(a(\sigma+\gamma)) = (d+1)/2.$  The result thus follows in
the case  $d\equiv3$ mod 4.  Next suppose that  $d\equiv1$ mod 4.
Then it may be checked directly
that $s(a)=u$, $s(a(\sigma +\gamma -1)) = u$, and
$s(a(\sigma+\gamma)) = 5(d+3)/12,$  and the result again follows.
\QED

\begin{Lemma} \label{GMWlem3}
Let $d=uv$ with $v=3$ and odd $u \geq 3$.  Then
every GMW difference set $G_{d,3,r}$ is inequivalent to the
difference set $D_{\sigma+\gamma, d}$ corresponding to the Glynn
type (I) hyperoval.
\end{Lemma}

{\Proof}
The result is easily checked for $u=3$ and $u=5$, because then the 2-ranks
of the two kinds of difference sets differ (see Table~1).
So let $u \geq 7$.
Write $q=2^u=2^{d/3}$.  Let $\chi$ be any nontrivial multiplicative
character on $\Ff_{2^d}^*$ whose restriction $\chi_1$ to
$\Ff_q^*$ is trivial.  Thus, if $\omega$ is a Teichm\"uller character
(of order $2^d -1$) on  $\Ff_{2^d}^*$ , then $\chi = \omega^{-a}$
for an integer $a$ divisible by $q-1$.

Suppose that $G_{d,3,r}$ were equivalent to $D_{\sigma+\gamma, d}$.
Then it would follow from Lemma~\ref{GMWlem1} and \eqref{eq2.1} that
$$J(\chi,\chi^{\sigma+\gamma -1}) = \mu g(\chi^t)$$
for some complex root of unity $\mu$ and some integer $t$ relatively
prime to $2^d -1$.  Setting $p=2$ and $d=3$ in the proof of
Theorem~\ref{thm1}, we see from that proof that $g(\chi^t)$
is divisible by $q$ in the ring of algebraic
integers.  Thus $J(\chi,\chi^{\sigma+\gamma -1})$ is divisible by
$q$,  so that every prime ideal occurring in the factorization of
$J(\chi,\chi^{\sigma+\gamma -1})$ over  $\Zz[\xi_{2^d-1}]$
has an exponent $\geq  u$.  Thus by \eqref{eq3.3a},
$$s(a)+s(a(\sigma +\gamma -1)) - s(a(\sigma+\gamma)) \geq u$$
for every integer $a$ divisible by $q-1$.  This contradicts
Lemma~\ref{GMWlem2}, and so the proof is complete.
\QED

\begin{Lemma} \label{conj4}
Let $u \geq 5$ be odd.  Then $A_6(3u)$ is never a power of\/ $3$
and $A_6(5u)$ is never a power of\/ $5$.
\end{Lemma}

{\Proof} By Theorem~\ref{thm:number-solutions}, for odd $d > 1$,
\begin{equation} \label{eq:A_6}
        A_6(d) = 2F_{(d-1)/2} - 1 .
\end{equation}
First suppose that $A_6(d)$ is a power of 3 for $d=3u$.  Then,
since $A_6(d)$ is different from $3,9,27$ (see Table~1),
$A_6(d)$ must be divisible by 81.
So by \eqref{eq:A_6}, $d=57+432x$ for some nonnegative integer $x$.  Since
$F_n$ (mod 109)
has period 108, it follows that $A_6(d)$ (mod 109) has period 216,
and hence $A_6(57+432x)$ is congruent to 42 modulo 109 for {\em all\/}
nonnegative
integers $x$.  This is a contradiction, since 42 is not a power of 3
modulo 109.

Next suppose that $A_6(d)$ is a power of 5 for $d=5u$.  Then 
since $A_6(d)$ is different from $5$ and $25$ (see Table~1), we have
$125\mid A_6(d)$.
So by \eqref{eq:A_6}, we have
$d=585+1000y$ for some nonnegative integer $y$.  Since $F_n$
(mod 251) has period 250, it follows that $A_6(d)$ (mod 251) has period
500, and so $A_6(585+1000y)$ is congruent to 235 modulo 251 for {\em all\/}
nonnegative $y$.  This is a contradiction, since 235 is not a power of 5
modulo 251.  
\QED

\begin{Theorem} \label{thm10}
The GMW difference set $G_{d,v,r}$ is inequivalent to each of
the difference sets  $D_{k,d}$ arising from the hyperovals in
{\em\eqref{reg-oval}}--{\em\eqref{Glynn-II-oval}}.
\end{Theorem}

{\Proof}
As was mentioned following the definition of $G_{d,v,r}$, the complement
$\overline{G_{d,v,r}}$ has 2-rank $uv^w$, where $d=uv$ and
$1<w<u<d$.
For $k$ as in \eqref{reg-oval} or \eqref{trans-oval}, the difference set
$\overline{D_{k,d}}$ has 2-rank $d$
(cf\@. the proof of Theorem~\ref{thm8}).  Since $u(d/u)^w > d$, it
follows that
$G_{d,v,r}$ is inequivalent to $D_{k,d}$.  

Now let $d$ be an odd composite
integer $> 5$.
It is easily checked by extending Table~1
that none of $A_6(d)$, $A_{\sigma+\gamma} (d)$, or $A_{3\sigma+4} (d)$
can equal
$(d/u)^{w-1}$ for any $d < 500$.  Hence assume $d > 500$.
{}From the proof of Lemma~\ref{lem:2-rank},
$$A_{3\sigma+4} (d) > \tfrac {1} {4}(2.09)^{(d+1)/2} > 3^{d/3} > (d/u)^{u-2} \ge
(d/u)^{w-1}.$$
This completes the proof for $k$ as in \eqref{Glynn-II-oval}.

Next, let $k$ be as in \eqref{Glynn-I-oval}.
By Lemma~\ref{GMWlem3},  $G_{d,3,r}$ cannot be equivalent to 
$D_{\sigma+\gamma, d}$.
Hence, assume that $u \le d/5$.  From the proof
of Lemma~\ref{lem:2-rank},
$$A_{\sigma+\gamma} (d) > \tfrac {1} {4}(1.92)^{(d+1)/2} > 5^{d/5} > (d/u)^{u-2}
 \ge
(d/u)^{w-1}.$$
This completes the proof for $k$ as in \eqref{Glynn-I-oval}.

Finally, let $k$ be as in \eqref{Segre-oval}, i.e., $k=6$.
By Lemma~\ref{conj4}, we cannot have $A_6(d)$ equal to $(d/u)^{w-1}$
with $u=d/3$ or $d/5$ or $d/9$.  Moreover, $A_6(d)$ cannot equal $(d/u)^{w-1}$
for $u=d/7$, because it follows from Table~1 and the recurrence
\eqref{eq:recurrence1} that
$A_6(d)$ never equals 7 and is never divisible by $49$ for any odd $d$.
Thus assume that $u \le d/11$.  From the proof of
Lemma~\ref{lem:2-rank},
$$A_6(d) > (1.6)^{(d+1)/2} >  11^{d/11} > (d/u)^{u-2} \ge  (d/u)^{w-1}.$$
This completes the proof for $k$ as in \eqref{Segre-oval}, and hence completes 
the proof of the theorem.
\QED

\noindent{\bf Remark:} 
In the preceding proof, in order to show the inequivalence of GMW
sets $G_{d,3,r}$ and the difference sets $D_{\sigma+\gamma,d}$
arising from the Glynn type (I) hyperovals, we used
Lemma~\ref{GMWlem3} instead of a 2-rank argument. We do in fact also
have a proof\footnote{It can be obtained on
request from the authors, or by WWW at {\tt
http://radon.mat.univie.ac.at/People/kratt/artikel/glynn.html}.}, 
which is similar in spirit to the 
proof of Lemma~\ref{conj4}, that for odd $u \geq 3$ the number
$A_{\sigma+\gamma} (3u)$ is never a power of\/ $3$, whence also the
2-ranks of $G_{d,3,r}$ and $D_{\sigma+\gamma,d}$ are distinct.
However, we prefer the approach via
Lemma~\ref{GMWlem3}, because 
the aforementioned proof requires extensive computer calculations,
and because 
Lemma~\ref{GMWlem1} is of independent interest.

\section{2-ranks of circulant matrices}

For integers $k \geq 3$ and $d \geq 2$, let 
$$f(x) = x^k +x^{k-1},\qquad x \in \Ff_{2^{d}}^*.$$
In Theorem~\ref{thm:gfunc}, we will determine the ranks over $\Ff_2$ of 
certain circulant matrices $M_k$, defined as follows.  Let
$\alpha$ denote a generator of the cyclic group $\Ff_{2^{d}}^*$
and define the circulant matrix
$$M_k = ( m_{i,j} )_{1 \leq i,j \leq 2^d - 1}\ ,$$
where $m_{i,j}$ denotes the number of $x \in \Ff_{2^{d}}^*$
for which $ \alpha^{j-i} = f(x)$.  
Note that for each $k$ such that
$D(x^k)$ is a hyperoval in $PG(2,2^d)$, we have $M_k = 2A_k$, where
$A_k$ is 
the incidence matrix (defined at the beginning of Section 3)
corresponding to the difference set $D_{k,d}$ (defined in
Theorem~\ref{thm2}), 
the 2-ranks of which we studied in Section~4.
This is because the incidence matrix $A_k$ is not affected by
changing $\tau(x)$ 
(as defined before Lemma~\ref{lem2})
from $x+x^k$ to $x^{k-1}+x^k$, since the difference
sets $D_{k,d}$ and $D_{k/(k-1),d}$ are equivalent by Lemma~\ref{lem3'}.

Let  ${\rm rank}_2(M_k)$ denote the rank of  $M_k$ over $\Ff_2$.
In Section~4 we were concerned with the computation of the 
2-ranks of incidence matrices $A_k$. Such a matrix $A_k$
can be viewed as an adjacency matrix of the directed graph where 
for two vertices $u$, $v$ in $\Ff_{2^d}^*$, a directed edge
connects $u$ to $v$ if and only if $v/u$ is in the difference 
set $D_{k,d}$, i.e., if and only if $v/u$ is in the image of $\tau$.
The following two theorems exhibit values of $k$
for which the computation of the 2-rank of $M_k$ is
also a 2-rank computation for an {\em adjacency} matrix $N_k$
of a directed graph, where this time $u$ is connected to $v$
if and only if $v/u$ is in the image of $f$.

\begin{Theorem}\label{thm:samerank1}
Let $d > 1$ be odd.  Let $N_5$ be the matrix obtained from $M_5$
by replacing every nonzero entry with $1$.  Then
${\rm rank}_2(N_5)$ =  ${\rm rank}_2(M_5)$.
\end{Theorem}

{\Proof}
It suffices to show that every nonzero entry of $M_5$ is odd.
We will prove the stronger result that for each 
$c \in \Ff_{2^{d}}^*$,  the polynomial $g(x) = x^5+x^4+c$ has 
either no zeros, one zero, or three zeros in $\Ff_{2^{d}}^*$.
First, if $g(x)$  had exactly two zeros in $\Ff_{2^{d}}^*$,
then $g(x)$ would have an irreducible cubic factor
(normal over  $\Ff_{2^{d}}$), so that $g(x)$ would have
five zeros in $\Ff_{2^{3d}}$ .  Similarly, if $g(x)$ had four zeros in 
$\Ff_{2^{d}}^*$, then $g(x)$ would have five zeros in $\Ff_{2^{d}}^*$.
Thus, since $3d$ is odd, it remains to show that $g(x)$ cannot have
fives zeros in  $\Ff_{2^{d}}$ for any odd $d$.
Assume for the purpose of contradiction that 
$g(x)$ has fives zeros in  $\Ff_{2^{d}}$.
Let $z$ denote one of these zeros.   Since $c$ is nonzero, $z$ 
is not 0 or 1.  Over  $\Ff_{2^{d}}$,
we have the factorization $g(x) = (x+z)h(x)$, where
$$h(x)=x^4+(z+1)x^3+(z^2+z)x^2+(z^3+z^2)x+(z^4+z^3).$$
Note that $h(0)$ is nonzero since $g(0)$ is nonzero.
We have
$$q(x) = (z^4+z^3)^{-1} x^4 h(z/x) = x^4+x^3+x^2+x+w,$$
where $w = z/(z+1)$.
Since $h(x)$ has four zeros in  $\Ff_{2^{d}}$, so does $q(x)$.
Write $q(x) = (x^2+rx+s)(x^2+ux+v)$ over  $\Ff_{2^{d}}$.
Since both of these quadratic factors are reducible over  $\Ff_{2^{d}}$,
we have $\Tr(s/r^2) = \Tr(v/u^2) = 0$.
Note that $r$ and $u$ are nonzero, since $q(x)$ cannot have multiple 
zeros over  $\Ff_{2^{d}}$.  We will obtain a contradiction by showing that
$\Tr(s/r^2)+\Tr(v/u^2) = 1$.   Now,
$r+u=1$, $ru+s+v=1$, and $su+vr = 1$.  Thus 
$r=1+u$, $u^2+u+v+1 = s$, and
$u^3+u^2+uv +u +(v+uv) = 1$.  Hence $v = 1+u+u^2+u^3$
and $s = u^3$.
Therefore $$s/r^2 = u^3 / (1+u^2) = u+1/(1+u)+1/(1+u)^2,$$
so that $\Tr(s/r^2) = \Tr(u)$.  On the other hand,
$v/u^2 = (1/u)+(1/u)^2 +1+u$, so that $\Tr(v/u^2) = 1+\Tr(u)$.
Thus $\Tr(v/u^2)+\Tr(s/r^2) = 1$, 
which yields the desired contradiction. This completes the proof.
\QED

The next theorem shows that the conclusion of Theorem~\ref{thm:samerank1}
is valid for {\em all} $d > 1$ when the subscript 5
is replaced by any $k$ of the form $2^m -1$ with $m > 1$.

\begin{Theorem}\label{thm:samerank2}
Let $d > 1$. Set $k = 2^m - 1$ for  an integer $m > 1$, and let $N_k$
be the matrix obtained from $M_k$ by replacing every nonzero entry with $1$.
Then ${\rm rank}_2(N_k)$ =  ${\rm rank}_2(M_k)$.
\end{Theorem}

{\Proof}
It suffices to show that
for each $c \in \Ff_{2^{d}}^*$, the polynomial $g(x) = x^k+x^{k-1}+c$
has either no zeros or an odd number of zeros in $\Ff_{2^{d}}^*$.
Assume that $z$ is a zero of $g(x)$ in  $\Ff_{2^{d}}^*$, 
and note that $z$ is not 0 or 1.
We will show that $g(x)$ has an odd number of zeros in  $\Ff_{2^{d}}^*$.
We have the factorization $g(x) = (x+z)h(x)$,
where  $$h(x) = x^{k-1}+(z+1)x^{k-2}+(z^2+z)x^{k-3}+\ldots +
(z^{k-2}+z^{k-3})x+(z^{k-1}+z^{k-2}).$$
It remains to show that $h(x)$ has an even number of zeros.  We have
$$q(x) = (z^{k-1}+z^{k-2})^{-1} x^{k-1} h(z/x) = 
x^{k-1}+x^{k-2}+\ldots+x+z/(z+1).$$
Since $k = 2^m -1$, it follows that
whenever $x$ is a zero of $q(x)$, so is $x+1$.  Thus $q(x)$ and
hence $h(x)$ has an even number of zeros. This completes the proof.
\QED

We next give a formula for  ${\rm rank}_2(M_k)$ in terms of $s(x)$,
where $s(x)$ is defined above Theorem~\ref{thm5}, with $q = 2^d$.  Define
$$s^*(x)=\begin{cases} s(x),&\text{if } 
(2^d-1) \hbox{${}\not\kern2.5pt\mid{}$} x,\\
d,&\text{if } (2^d-1) \mid x.\end{cases}$$

\begin{Theorem}\label{thm:rankformula}
For integers $k \geq 3$ and $d \geq 2$,  ${\rm rank}_2(M_k)$ equals
the number of $a$, $0 < a < 2^d-1$, for which
\begin{equation}\label{eq:eqs^*}
s(a)+s((k-1)a) = s^*(ka).
\end{equation}
\end{Theorem}

{\Proof}
The $2^d -1$ complex eigenvalues of $M_k$ are
$$\sum_{x\in \Ff_{2^{d}}^*}\chi(f(x)),$$
where $\chi$ ranges over the $2^d -1$ multiplicative characters
on $\Ff_{2^{d}}^*$.  Let 
$\mathfrak{p}$ be a prime ideal in $\Zz[\xi_{2^d-1}]$ lying over
$2$, and let $\omega$ =
$\omega_{\mathfrak p}$ be the Teichm\"uller character on $\Ff_{2^d}^*$.
Then  ${\rm rank}_2(M_k)$ equals
the number of $a$, $0 < a < 2^d-1$, for which
$$\mathfrak{p} \hbox{${}\not\kern2.5pt\mid{}$} 
J(\omega^{-a},\omega^{-(k-1)a}) =
\sum_{x\in \Ff_{2^{d}}^*}\omega^{-a}(f(x))$$
(cf. Theorem~\ref{thm4}).
For $0 < a < 2^d-1$, we have $$
\mathfrak{p}^{s(a)+s((k-1)a)-s^*(ka)}\dmid J(\omega^{-a},\omega^{-(k-1)a})
$$ (cf. \eqref{eq3.3a}),  since
$J(\omega^{-a},\omega^{-(k-1)a})$ is a unit if $ka$ or $(k-1)a$ is
divisible by $2^d - 1$.
The theorem now follows.
\QED

Let  $R_k(d)$ denote the number of solutions $a$, $0 < a < 2^d-1$,
to equation \eqref{eq:eqs^*}.  Thus by Theorem~\ref{thm:rankformula}, 
$R_k(d) = {\rm rank}_2(M_k)$.  The reader should observe the
similarity of the problems of determining $R_k(d)$,
the number of solutions $a$, $0 < a < 2^d-1$,
to equation \eqref{eq:eqs^*}, and of determining $B_k(d)$,
the number of solutions $a$, $0 < a < 2^d-1$,
to equation \eqref{eq:keyeq1}. In the case of $B_k(d)=dA_k(d)$ we obtained 
linear
recurrences for all relevant values of $k$ (see
\eqref{eq:recurrence1}, \eqref{rec:Glynn1}, \eqref{rec:Glynn2}).
The next proposition shows that linear recurrences exist for the
$R_k(d)$ as well.  We omit the proof, as the proposition can be 
proved by an approach similar to the one we used in the proofs of
Theorems~\ref{conj3} and \ref{conj2}. 

\begin{Proposition} \label{prop:rec}
Let $k$ be a fixed integer $\ge3$.
The generating function $\sum_{d \geq 2} R_k(d) z^d$ for the numbers
$R_k(d)$ of solutions $a$, $0 < a < 2^d-1$, to equation 
{\em\eqref{eq:eqs^*}} is a
rational function in $z$, and, hence, the sequence $(R_k(d))_{d \ge e}$
satisfies a linear recurrence for $e$ large enough.
\QED
\end{Proposition}
A proof of this proposition in the style of the proofs of
Theorems~\ref{conj3} and \ref{conj2} also gives bounds on
the degrees of the numerator and denominator polynomials of the
rational generating function $\sum_{d \geq 2} R_k(d) z^d$.
So, this allows us, at least for $k$ not too large, to determine
recurrences via the computer, by first computing enough numbers $R_k(d)$
to guess a linear recurrence using {\tt gfun} \cite{SaZi}, and then
applying Lemma~\ref{lem:rec-always}.

In the next theorem, we determine $R_k(d)$
explicitly for $3 \leq k \leq 9$ by computing the corresponding
rational generating functions.

\begin{Theorem} \label{thm:gfunc}
For $3 \leq k \leq 9$, the generating functions $\sum_{d \geq 2} 
R_k(d) z^d$ are respectively
\begin{align} \label{eq:GF3}
\sum_{d \geq 2} R_3(d) z^d&=\frac {z^2 (2-z)} 
{(1-z)(1-z-z^2)},\\
\label{eq:GF4}
\sum_{d \geq 2} R_4(d) z^d&=\frac {2z^2} 
{1-z^2},\\
\label{eq:GF5}
\sum_{d \geq 2} R_5(d) z^d&=\frac {z^3 (3+2z-3z^2)} 
{(1-z)(1+z^2)(1-z-z^2)},\\
\label{eq:GF6}
\sum_{d \geq 2} R_6(d) z^d&=\frac {z^2 (2+4z^2)} 
{1-z^2-z^4},\\
\label{eq:GF7}
\sum_{d \geq 2} R_7(d) z^d&=\frac {z^2 (2+2z-6z^2-2z^3-2z^4+5z^5)} 
{(1-z)(1-z-z^3)(1-z^2-z^3)},\\
\label{eq:GF8}
\sum_{d \geq 2} R_8(d) z^d&=\frac {6z^3} 
{1-z^3},\\
\label{eq:GF9}
\sum_{d \geq 2} R_9(d) z^d&=\frac {z^2(2-4 z+6 z^2 +2 z^3 +2 z^4 -12 z^5 -
2z^6 +7 z^7 )} 
{(1-z) (1-z-z^2) (1+z^3-z^6)}.
\end{align}

\end{Theorem}

{\Proof} 
These assertions can, in principle, be proved as we described above the
statement of the theorem.
However, it is 
instructive to exhibit the solutions to \eqref{eq:eqs^*} in each case
explicitly, from which \eqref{eq:GF3}--\eqref{eq:GF9} then follow
easily by means of generating function calculus.

For the description of the solution sets we use a standard notation
from the calculus of words (cf\@. \cite{loth}): 
Given an alphabet $\mathcal A$, the set of all words 
(including the empty word) consisting of
letters from $\mathcal A$ is denoted by $\mathcal A^*$.
With that terminology, the solutions to \eqref{eq:eqs^*} in
the cases $k=3,\dots,9$ can be described as follows:

For $k=3$, the binary representations of the 
solutions $a$ to \eqref{eq:eqs^*}, $0<a<2^d-1$, 
are obtained by forming all possible rotations of the strings 
of length $d$
from the set $\{0,01\}^*\backslash\{0\}^*$.

For $k=4$, the binary representations of the 
solutions $a$ to \eqref{eq:eqs^*}, $0<a<2^d-1$,  
are obtained by forming all possible rotations of the strings 
of length $d$
from the set $\{01\}^*$.

For $k=5$, the binary representations of the 
solutions $a$ to \eqref{eq:eqs^*}, $0<a<2^d-1$,  
are obtained by forming all possible rotations of the strings 
of length $d$
from the set $\{0,001,0011\}^*\backslash\{0\}^*$.

For $k=6$, the binary representations of the 
solutions $a$ to \eqref{eq:eqs^*}, $0<a<2^d-1$,  
are obtained by forming all possible rotations of the strings 
of length $d$
from the set $\{01,0011\}^*$.

For $k=7$, the binary representations of the 
solutions $a$ to \eqref{eq:eqs^*}, $0<a<2^d-1$,
are obtained by forming all possible rotations of the strings 
of length $d$
from the set $(\{0,001\}^*\cup\{01,011\}^*)\backslash\{0\}^*$.

For $k=8$, the binary representations of the 
solutions $a$ to \eqref{eq:eqs^*}, $0<a<2^d-1$,
are obtained by forming all possible rotations of the strings 
of length $d$
from the set $\{001,011\}^*$.

For $k=9$, the binary representations of the 
solutions $a$ to \eqref{eq:eqs^*}, $0<a<2^d-1$,  
are obtained by forming all possible rotations of the strings 
of length $d$
from the set $(\{0,00011,000111\}\cup\{0001,000101,00010101,\dots\})^*
\backslash\{0\}^*$.

These assertions can be established in a way similar to the proof of
Proposition~\ref{prop:solutions}, with use of the fact that the 
solutions   $a$ to \eqref{eq:eqs^*}, $0<a<2^d-1$, are characterized
by the property that there is no instance of a 1 occurring in the
same place in the binary expansions of $a$ and $(k-1)a$ (mod $2^d-1$).

The formulas \eqref{eq:GF3}--\eqref{eq:GF9} for the generating
functions follow, as we now demonstrate for the case $k=5$.

{}From the above, the binary representation of a
solution $a$ to \eqref{eq:eqs^*} with $k=5$ must be some rotation of a
nonzero string of length $d$ obtained by concatenating the blocks $0$, $001$,
and $0011$.  
To identify the block in the {\em cyclic} binary representation of $a$
which overlaps with the rightmost place of the binary representation of
$a$, we underline the digit in that block where the overlap occurs.
We have the following possibilities for such blocks:
$\underline0$,
$\underline001$,
$0\underline01$,
$00\underline1$,
$\underline0011$,
$0\underline011$,
$00\underline11$,
$001\underline1$.
Once we have chosen one of these possibilities, we have to ``glue" this
already chosen block with a (possibly empty) sequence of other blocks from
$\{0,001,0011\}$ to form a  binary string representing a
solution to \eqref{eq:eqs^*}, where the underlined place is considered
as the rightmost place in the string (through appropriate rotation).
Hence, in symbolic notation again, 
the set of all solutions to \eqref{eq:eqs^*} is
$$\big(\{\underline0,
\underline001,
0\underline01,
00\underline1,
\underline0011,
0\underline011,
00\underline11,
001\underline1\}\times\{0,001,0011\}^*\big)
\backslash\{0\}^*.$$
In particular,
a string taken from this set is a solution
to \eqref{eq:eqs^*} with corresponding $d$ equal to the length of
the string.
Since the set 
$\{\underline0,
\underline001,
0\underline01,
00\underline1,
\underline0011,
0\underline011,
00\underline11,
001\underline1\}$
contains one string of length 1, three strings of length 3, and four
strings of length 4, its corresponding generating function is
$z+3z^3+4z^4.$  Similarly, the set $\{0,001,0011\}$ has generating
function $z+z^3+z^4.$
Therefore, by elementary
``generatingfunctionology" (see \cite{wilf}), 
the generating function for 
the solutions $a$ to
\eqref{eq:eqs^*}, and, hence, the generating function 
$\sum_{d \geq 2} R_5(d) z^d$, is given by
\begin{equation} \label{eq:GF5a}
(z+3z^3+4z^4)\cdot
\frac {1} {1-(z+z^3+z^4)}-\frac {z} {1-z}.
\end{equation}
Simplification of expression \eqref{eq:GF5a} gives \eqref{eq:GF5}.

We leave the details in the remaining cases to the reader.
\QED

\bigskip
{\bf Acknowledgements}.
The authors are grateful to S. Payne for pointing out 
\cite[Theorem~12.5.3]{pt} in the proof of Lemma~\ref{lem3'}.
The second and fourth authors would like to express their gratitude to the
Mathematische Forschungsinstitut Oberwolfach, where they first met and where
part of this work initiated. 
The third author would like to thank the
University of California at San Diego for the opportunity to spend a
visit in Winter Quarter 1998, which brought together the first and
third authors to work on this subject.  The research of the fourth
was supported in part by NSA grant MDA 904-97-0104, and a grant from the 
University of Delaware Research Foundation.

\bigskip


\begin{thebibliography}{99}

\bibitem{ab}M. Antweiler and L. B\"omer, {\em Complex sequences
over $GF(p^m)$ with a two-level autocorrelation function and a
large linear span}, IEEE Trans\@. Inform\@. Theory, {\bf 38}
(1992), 120--130. 

\bibitem{asskey}
E. F. Assmus, Jr., and J. D. Key, {\em  Designs and their codes},
Cambridge Tracts in
Mathematics, 103, Cambridge University Press, Cambridge, 1992.

\bibitem{be} B. C. Berndt and R. J. Evans,
{\it Sums of Gauss, Jacobi, and Jacobsthal\/},
J. Number Theory, {\bf 11} (1979), 349--398.

\bibitem{be1} B. C. Berndt, R. J. Evans and K. S. Williams,
{\em Gauss and Jacobi sums}, Wiley Interscience, 1998.


\bibitem{BrHaHa} W. G. Bridges, M. Hall, and J. L. Hayden, {\em Codes and
designs}, J. Combin\@. Theory Ser.~A
 {\bf 31} (1981), 155--174.

\bibitem{cp} J. Cannon and C. Playoust, {\em An introduction to MAGMA},
University of Sydney, Australia, 1993.  


\bibitem{cs} W. E. Cherowitzo and L. Storme, {\it $\alpha$-flocks with
oval herds and monomial hyperovals\/}, preprint.

\bibitem{gl}D. Glynn, {\it Two new sequences of ovals in finite
Desarguesian planes of even order\/}, Lecture Notes in Mathematics
1036, Springer-Verlag, 1983, 217--229. 

\bibitem{gd} J. M. Goethals and P. Delsarte, {\it On a class of
majority logic decodable cyclic codes\/}, IEEE Trans\@. Inform\@. Theory,
{\bf 14} (1968), 182--188.  

\bibitem{gol} S. W. Golomb, {\it The use of combinatorial structures
in communication signal designs\/}, in: ``Applications of Combinatorial
Mathematics'', Chris Mitchell (ed.), vol.~60, 1997, 59--78.  

\bibitem{go} B. Gordon, W. H. Mills and L. R. Welch, {\it Some new
difference sets\/}, Canad\@. J. Math\@. {\bf 14} (1962), 614--625.   



\bibitem{h} N. Hamada, {\it On the $p$-rank of the incidence matrix of
a balanced or partially balanced incomplete block design and its
applications to error-correcting codes\/}, Hiroshima Math. J. {\bf 3}
(1973), 154--226. 

\bibitem{ho} N. Hamada and H. Ohmori, {\it On the BIB-design having
the minimum $p$-rank\/}, J. Combin\@. Theory Ser.~A {\bf 18} (1975),
131--140. 

\bibitem{hi} J. W. P. Hirschfeld, {\it Projective Geometries over
Finite Fields}. Oxford University Press, 1979. 

\bibitem{ju} D. Jungnickel, {\it Difference Sets\/}, in: 
``Contemporary Design Theory, A Collection of
Surveys'', J.~Dinitz, D.~R.~Stinson (eds.),  Wiley-Interscience Series
in Discrete Mathematics and 
Optimization, Wiley, New York, 1992, 241--324.
  
\bibitem{lan} E. S. Lander, {\it Symmetric Designs, An Algebraic
Approach\/}, Cambridge University Press, Cambridge, 1983. 

\bibitem{lidl} R. Lidl and H. Niederreiter, {\it Finite Fields}, 2nd~ed.,
Encyclopedia of Mathematics and Its Applications, vol.~20, Cambridge
University Press, Cambridge, 1997.
  
\bibitem{loth} M. Lothaire, {\it Combinatorics on words},
2nd~ed., Encyclopedia of
Mathematics and Its Applications, vol.~17, Cambridge
University Press, Cambridge, 1997.

\bibitem{mm} J. MacWilliams and H. B. Mann, {\it On the $p$-rank of
the design matrix of a difference set\/}, Inform\@. Control {\bf 12}
(1968), 474--488. 

\bibitem{MacWiSl} F. J. MacWilliams and N. J. A. Sloane, {\it The theory of
error-correcting codes}, North-Holland, 1977. 

\bibitem{ma} A. Maschietti,
{\it Difference sets and hyperovals\/},
Designs, Codes and Crypt\@. {\bf 14} (1998), 89--98.

\bibitem{pt} S. E. Payne and J. A. Thas, {\it Finite Generalized Quadrangles\/}, 
Research Notes in Mathematics, {\bf 110}, Pitman Pub\@. Inc\@., 1984.

\bibitem{SaZi}
B. Salvy and P. Zimmermann, {\em Gfun: a Maple package for the
manipulation of generating and holonomic functions in one 
       variable}, ACM Trans\@. Math\@.
       Software {\bf 20} (1994).

\bibitem{schwel}R. A. Scholtz and L. R. Welch, {\it GMW sequences}, 
IEEE Trans\@. Inform\@. Theory {\bf 30} (1984), 548--553.

\bibitem{segre}
B. Segre, {\em Ovals in a finite projective plane}, 
Canad\@. J. Math\@. {\bf 7} (1955),
414--416.

\bibitem{sm} K. J. C. Smith, {\it On the $p$-rank of the incidence
matrix of points and hyperplanes in a finite projective geometry\/},
J. Combin\@. Theory {\bf 7} (1969), 122--129. 

\bibitem{Stanley} R. P. Stanley, {\it Enumerative Combinatorics},
Vol.~1, Wadsworth \& Brooks/Cole, Pacific Grove, California, 1986;
reprinted by Cambridge University Press, Cambridge, 1997.

\bibitem{turyn}R. J. Turyn, {\it Character sums and difference sets}, 
Pacific J. Math\@. {\bf 15} (1965), 319--346.


\bibitem{wa} L. C. Washington, {\it Introduction to Cyclotomic
Fields\/}, 2nd edition, Springer, 1997. 
  
\bibitem{wilf} H. S. Wilf, {\em generatingfunctionology}, 2nd ed.,
Academic Press, San Diego, 1994.

\bibitem{qx} Q. Xiang, {\it On balanced binary sequences with
two-level autocorrelation functions\/}, 
IEEE Trans\@. Inform\@. Theory, to appear.  


\bibitem{kYa1} K. Yamamoto, 
{\it On Jacobi sums and difference sets\/}, J. Combin\@. Theory Ser.~A
{\bf 3} (1967), 146--181. 


\bibitem{kYa2} K. Yamamoto, {\it On congruences arising from relative
Gauss sums\/}, 
in: ``Number Theory and Combinatorics'', Japan 1984, World Scientific
Publ., 1985, 423--446. 


\end{thebibliography}
\end{document}